\documentclass[10pt]{article}
\usepackage{graphicx,epsfig}
\usepackage{amsmath}

\usepackage{amsmath}
\usepackage[makeroom]{cancel}
 \usepackage{color}      

\usepackage{url}

\usepackage{tabu}

\usepackage{array}
\usepackage{makecell}

\usepackage{multirow}

\definecolor{auburn}{rgb}{0.43, 0.21, 0.1}

\def\dyad{\! \otimes \!}

\usepackage{subcaption}
\usepackage[driverfallback=dvipdfm,
bookmarks=true,    
bookmarksopen=true,
bookmarksopenlevel=\maxdimen,
bookmarksdepth=10
]{hyperref}
\hypersetup{
pdfcenterwindow=true, 
pdfpagelabels=true,
pagebackref=true,            
hypertexnames=true,
plainpages=false,
unicode=false,     
pdftoolbar=true,                
    pdfmenubar=true,          
    pdffitwindow=false,         
    pdfstartview={FitH},        
    pdftitle={On Third-Order Accuracy of MUSCL Scheme},
    pdfauthor={Hiroaki Nishikawa},    
    pdfsubject={On Third-Order Accuracy of MUSCL Scheme},       
    pdfcreator={Hiroaki Nishikawa},   
    pdfproducer={Hiroaki Nishikawa}, 
    baseurl={http://www.hiroakinishikawa.com},
    pdfkeywords={}, 
    pdfnewwindow=true,     
    colorlinks=true,       
    linkcolor=black,       
    citecolor=black,       
    filecolor=black,       
    urlcolor=black,        
  breaklinks=true,
  hyperfigures=true,
  backref=true,
breaklinks=false, 
pdfpagelabels,      
pagebackref,        
hypertexnames=true, 
plainpages=false,   
naturalnames        
}
\usepackage[all]{hypcap}

\usepackage[thinlines]{easytable}

\usepackage{siunitx}

\usepackage{color}      
\usepackage{placeins}
\usepackage{indentfirst}
\usepackage{multirow}

\usepackage{amssymb}
\usepackage{amsmath,latexsym}
\usepackage{url}
\usepackage{enumerate}
\usepackage{graphicx}

\usepackage{booktabs}
 \newcommand{\ra}[1]{\renewcommand{\arraystretch}{#1}}
 
\usepackage{subcaption} 
\usepackage{mathtools}

\usepackage{cleveref}

\def\dyad{\! \otimes \!}

\usepackage{cleveref}
\usepackage[thinlines]{easytable}

\usepackage{blindtext,graphicx}
\usepackage[absolute]{textpos}
\definecolor{black}{rgb}{0.43, 0.21, 0.1} 
\usepackage{enumitem,xcolor}

\begin{document}

\title{\bf {\color{black} An Efficient Quadratic Interpolation Scheme for a Third-Order Cell-Centered Finite-Volume Method on Tetrahedral Grids} }

 \author{
{Hiroaki Nishikawa}\\
  {\itshape 
{National Institute of Aerospace}, Hampton, VA 23666, USA}
 \\
\and {Jeffery A. White}
  \\
\centerline{\itshape NASA Langley Research Center, Hampton VA 23681}%
}

\date{\today}
\maketitle

\begin{abstract}
In this paper, we {\color{black} propose an efficient quadratic interpolation formula utilizing solution gradients computed and stored at nodes and demonstrate its application to a third-order cell-centered finite-volume discretization on tetrahedral grids. The proposed quadratic formula is constructed based on an efficient formula of computing a projected derivative. It is efficient in that it completely eliminates the need to compute and store second derivatives of solution variables or any other quantities, which are typically required in upgrading a second-order cell-centered unstructured-grid finite-volume discretization to third-order accuracy. Moreover, a high-order flux quadrature formula, as required for third-order accuracy, can also be simplified by utilizing the efficient projected-derivative formula, resulting in a numerical flux at a face centroid plus a curvature correction not involving second derivatives of the flux. Similarly, a source term can be integrated over a cell to high-order in the form of the source term evaluated at the cell centroid plus a curvature correction, again, not requiring second derivatives of the source term. The discretization is defined as an approximation to an integral form of a conservation law but the numerical solution is defined as a point value at a cell center, leading to another feature that there is no need to compute and store geometric moments for a quadratic polynomial to preserve a cell average.} Third-order accuracy and improved second-order accuracy are demonstrated and investigated for simple but illustrative test cases in three dimensions. 
\end{abstract}

\section{Introduction}

Many compressible-flow problems are convection-dominated and accurate numerical simulations of such problems require superior accuracy in convection schemes. {\color{black} Modern high-order methods such as continuous and discontinuous Galerkin methods can certainly provide high-order accuracy and thus may be employed if a major code restructuring is allowed and high-order meshes are available for target problems. On the other hand, for the sake of simplicity in implementation and application with existing linear meshes}, there have also been continuing efforts in developing economical, low-dissipation, and third-order schemes for practical unstructured-grid computational fluid dynamics (CFD) codes \cite{AbalakinKozubskayaDervieux:IJA2014,yang_harris:AIAAJ2016,oswald_etal:IJNMF2016,liu_nishikawa_aiaa2016-3969,nishikawa_liu_aiaa2018-4166}. However, for second-order cell-centered finite-volume methods widely used in practical CFD codes, an extension to third-order accuracy is not a simple task. Third-order reconstruction-based cell-centered finite-volume methods require various algorithmic modifications: (i) computing cell-averages of source terms instead of a point evaluation, (ii) implementing a high-order flux quadrature formula with multiple numerical flux evaluations per face, (iii) a quadratic least-squares (LSQ) method involving more neighbor cells, (iv) computing and storing geometric moments in each cell for a quadratic solution reconstruction/interpolation, (v) computing second derivatives of all solution variables at every residual calculation, and (vi) computing cell-averaged exact solutions to perform code verification (e.g., see Refs.~\cite{barth_frederickson_AIAA1990,DelanayeEssers:AIAAJ1997,GoochVanAltena:JCP2002,CaraeniHill_AIAAJ2010,HaiderBrennerCourbetCroisille2011,charest_groth_etal:ICCFD2012,jalali_gooch:CF2017,PintBrennerCinnellaMaugarsRobinet:JCP2017,SetzweinEssGerlinger:JCP2021}). {\color{black}  If one wishes to develop an arbitrarily high-order solver beyond third-order accuracy, it might be worthwhile performing all these modifications because these would make the extension to fourth- and higher-order accuracy relatively easier \cite{jalali_gooch:CF2017}. However, it would be also highly desirable to reduce the burden to upgrade an existing code to third-order accuracy if one wishes to quickly improve an existing second-order CFD code.} In this paper, we will show that the extension can actually be made simpler and some of these modifications, including computations and storage for second derivatives, are not required. We also show that the accuracy of a second-order finite-volume code can be improved relatively easily by an economical quadratic interpolation scheme at almost no additional cost. The simplification is made possible by {\color{black} a key idea of eliminating second derivatives from curvature terms and the objective of this paper is to introduce this idea.} 

Traditionally, cell-centered finite-volume discretizations are constructed based on the integral form of a target conservation law: 
\begin{eqnarray}
\frac{ \partial \overline{\bf u}_j }{\partial t}  +  \frac{1}{V_j} \oint_{\partial V_j} {\bf f} \, ds    = \overline{\bf s}_j, 
\label{int_form00}
\end{eqnarray}
where $t$ is time, $V_j$ is the measure of a control volume $j$,  $\partial V_{j}$ denotes the boundary of the control volume $j$, $ds$ is an infinitesimal area of the boundary, $\overline{\bf u}_j$ is a vector of cell-averaged conservative variables, $\overline{\bf s}_j$ is a cell-averaged source term, and ${\bf f}$ is the flux projected along the control-volume face normal direction. The integral form is an exact relationship and therefore, one can construct an arbitrarily high-order discretization by discretizing the flux integral and the source term to the required order of accuracy with cell-averaged solutions stored at cells. However, a second-order accurate discretization 
is {\color{black} usually implemented for the sake of simplicity and cost reduction}, with point-valued solutions and source terms as,
\begin{eqnarray}
\frac{ d {\bf u}_j }{d t}   +  \frac{1}{V_j} \sum_{k \in \{  k_j \}}   {\Phi}_{jk}  = {\bf s}_j,
\label{2nd-order_CCFV}
\end{eqnarray}
where ${\bf u}_j$ and ${\bf s}_j$ are point values at a cell center ${\bf x}_j$, $\{  k_j \}$ is a set of faces of the cell $j$, and ${\Phi}_{jk}$ is a numerical flux evaluated at a face center. Treating the cell-averages as point values still leads to second-order accuracy as demonstrated in Ref.~\cite{nishikawa_centroid:JCP2020}, but it will not lead to third- or higher-order accuracy. Therefore, to extend it to third-order accuracy in a traditional manner, one will have to (i) treat ${\bf u}_j$ precisely as the cell average $\overline{\bf u}_j$, (ii) implement a high-order volume quadrature formula for computing $\overline{\bf s}_j$ (replacing ${\bf s}_j$), (iii) replace the single flux evaluation by a high-order surface quadrature formula requiring multiple evaluations of the numerical flux, (iv) implement a quadratic solution reconstruction scheme, (v) implement a quadratic LSQ gradient method, and (vi) compute and store second derivatives of the solution variables. Moreover, the quadratic reconstruction will require the computation and storage for geometric moments: e.g., monomials integrated over a cell as shown below in a quadratic polynomial of a scalar variable $u$, 
\begin{eqnarray}
u ({\bf x})  =  \overline{u}_j  
+     \left[ (  {\bf x} - {\bf x}_j ) \cdot  \nabla -    \int_{V_j}  \left\{  (  {\bf x} - {\bf x}_j ) \cdot  \nabla   \right\}  \, dV  \right]      u_j  
+     \frac{1}{2}   \left[   \left\{ ( {\bf x} - {\bf x}_j ) \cdot  \nabla  \right\}^2   
                                     -  \int_{V_j}  \left\{ ( {\bf x} - {\bf x}_j ) \cdot  \nabla  \right\}^2     dV       \right]    u_j ,
                                     \label{quadratic_MUSCL_reconst}
\end{eqnarray}
 such that the cell average of the quadratic polynomial reduces to the cell averaged solution $\overline{u}_j$. While the volume integral of $ ({\bf x} - {\bf x}_j ) \cdot  \nabla $ can be made to vanish by choosing ${\bf x}_j$ to be the geometric centroid, the volume integral of $ \left\{ ( {\bf x} - {\bf x}_j ) \cdot  \nabla  \right\}^2$ will not vanish in general and thus will generate geometric moments to be computed and stored, e.g., $ \int_{V_j}  ( {x} - {x}_j )^2 dV $. In addition, one would have to implement a mechanism to compute cell averages of exact solutions as well as source terms in the method of manufactured solutions for code verification purposes. Therefore, it requires a substantial effort to upgrade a second-order finite-volume code to third-order accuracy.

{\color{black}
To simplify the discretization, we first define the numerical solutions as point values at cell centers rather than cell averages. This type of finite-volume discretization corresponds to the QUICK (Quadratic Upstream Interpolation for Convective Kinematics) scheme of Leonard \cite{Leonard_QUICK_CMAME1979}, which is also known as a deconvolution finite-volume method \cite{FeliceDenaroMoela:NHT1993,Denaro:IJNMF1996,Denaro_CMAME:2015}. In this approach, we will need to compute a high-order accurate point-valued solution at a face from point-valued solutions stored at cells, not from cell-averaged solutions. Therefore, we need a quadratic interpolation rather than the quadratic solution reconstruction (\ref{quadratic_MUSCL_reconst}), and therefore, the geometric terms are not needed, saving computational time and memory. See Ref.~\cite{Nishikawa_3rdQUICK:2020} for fundamental differences between this approach and a finite-volume method with cell averages. To achieve third-order accuracy, we will still have to integrate the flux over a face with a high-order quadrature formula, and also apply another high-order quadrature formula to the source terms including the time derivative terms. As we will show, all these quadratures can be performed in the form of a face/cell center value plus a curvature correction, not requiring multiple quadrature points. This curvature correction approach will simplify the implementation because the face/cell center values already exist in a second-order code, thereby only requiring the addition of curvature corrections. These correction terms require second derivatives of the fluxes, solutions, and source terms, but as we will show in this paper, these second derivatives can be expressed in terms of first derivatives, thereby completely eliminating the need to compute and store second derivatives. The elimination of these second-derivative terms will lead to significant savings in computing time and memory, and also greatly simplify the implementation of the third-order discretization.  

The key idea is to extend the face-averaged nodal-gradient (FANG) cell-centered finite-volume method \cite{NishikawaWhite_FANG:jcp2020,Nishikawa_FANG_AQ:Aviation2020,WhiteNishikawaBaurle_scitech2020,NishikawaWhite:AIAA2021-2720} in a manner that allows the construction of a quadratic interpolation scheme that does not require second derivatives. In the FANG approach, we compute and store solution gradients at nodes from solution values stored at cells, and then perform a solution interpolation at a face by using the gradient averaged over the face nodes. This method has been demonstrated to reduce the cost of gradient computation and memory requirement for simplex-element grids, and reduce the residual stencil size \cite{NishikawaWhite_FANG:jcp2020,WhiteNishikawaBaurle_scitech2020}. Its robustness has been demonstrated for hypersonic-flow applications \cite{WhiteNishikawaBaurle_scitech2020,VULCAN_FANGplus_scitech2022}. As we will show, with gradients available at nodes, second derivatives required in a quadratic interpolation scheme can be evaluated using \color{black} the Green-Gauss} gradient scheme \cite{Frink:AIAAJ1998,pandya_etal:AIAAJ2016} over a cell and expressed entirely in terms of the nodal gradients, thereby eliminating the need to explicitly compute and store second derivatives. This third-order accurate scheme is simple to implement and computationally inexpensive compared with conventional quadratic reconstruction schemes because the total number of second derivatives to be computed and stored will add up to 30 in each cell for the Euler and Navier-Stokes equations having five variables in three dimensions and there is no need to compute and store them as well as any geometric moments. The quadratic interpolation without second derivatives is achieved by a generalization of the technique that has been utilized in the NASA USM3D code  \cite{Frink:AIAAJ1998,pandya_etal:AIAAJ2016}, where a linear solution reconstruction is performed efficiently in a tetrahedral cell by first recovering solution values at nodes from those at the cells and using these nodal solution values in the divergence theorem to compute the cell Green-Gauss solution gradients. The resulting expression can be expressed entirely by the nodal solutions and thereby does not require the actual computation of the cell gradient. In this paper, we will apply this idea to the quadratic term in a high-order interpolation scheme. This idea can be generalized to an arbitrary cell and an arbitrary interpolation point, however, in this work, we will focus on tetrahedral grids, due to their simplicity as well as our interest in anisotropic grid adaptation, which is known to be most efficiently performed with tetrahedral grids \cite{park_darmofal:AIAA2008-917,AlauzetDervieuxLoseille:JCP2010,Kleb_etal_aiaa2019-2948}. 

It should be noted that there exists a similar economical third-order discretization method for a node-centered formulation: the third-order edge-based discretization method \cite{liu_nishikawa_aiaa2016-3969,Katz_Sankaran_JCP:2011,katz_sankaran:JSC_DOI,diskin_thomas:AIAA2012-0609,LiuNishikawa_2017-3443}, which also does not require second derivatives nor any geometric moments. However, it achieves third-order accuracy only on simplex-element grids, i.e., triangular and tetrahedral grids in two and three dimensions, respectively. {\color{black} The third-order edge-based method has recently been extended to a cell-centered method \cite{KongDongLiu_EB3:jcp2023}, but again third-order accuracy is achieved only on simplex-element grids. On the other hand, the method proposed in this paper has the potential to achieve third-order accuracy on grids of arbitrary cell types.} 



This paper will focus on the key idea of the proposed efficient quadratic interpolation formula and its application to a third-order cell-centered finite-volume method, presenting numerical results for relatively simple but illustrative examples. Further algorithmic developments (e.g., nonoscillatory versions and extensions to grids of other cell types) and applications to practical low- and high-speed turbulent-flow problems will be addressed in future work. 

In Section \ref{sec_Conservative_Differential_Form},  we will discuss {\color{black} the baseline finite-volume discretization with point-valued solutions}.
In Section \ref{Efficient_Projected_Derivative_Formula_via_Nodal_Values}, we recall a well-known algorithm used for efficiently performing linear reconstruction over a tetrahedral cell and consider its generalization, which will be used repeatedly to eliminate second derivatives from a third-order accurate discretization. 
In Section \ref{Economical_Third_Order_Cell_Centered_Discretization}, we will describe the proposed third-order accurate cell-centered finite-volume discretization, extending the FANG cell-centered finite-volume method. {\color{black} In Section \ref{Solver_TimeScheme}, we will describe a nonlinear solver used for steady problems, a time integration scheme used for an unsteady problem, and a fixed-point iteration scheme used to approximately invert the mass matrix arising from a high-order quadrature applied to the time derivative terms.}
In Section \ref{Results}, numerical results are shown. 
Finally, in Section \ref{conclusions}, the paper concludes with remarks.


\section{{\color{black} Baseline Cell-Centered Finite-Volume Discretization with Point-Valued Solutions} }
\label{sec_Conservative_Differential_Form}
\indent 
\indent 
Consider a conservation law:
\begin{eqnarray}
\frac{ \partial {\bf u} }{\partial t}  + \mbox{div}{\cal F} = {\bf s}, 
\label{diff_form}
\end{eqnarray}
where ${\cal F}$ is a flux tensor.  In this paper, we
focus on the Euler equations: 
\begin{eqnarray}
{\cal F}  =   \left[\begin{array}{c}
\displaystyle  \rho {\bf v}  \\  [1.7ex]
\displaystyle  \rho  {\bf v} \dyad  {\bf v} +  p {\bf I}  \\   [1.7ex]
\displaystyle  \rho {\bf v} H 
\end{array}\right], 
\quad
\end{eqnarray}
where $\dyad$ denotes the dyadic product, ${\bf I}$ is the 3$\times$3 identity matrix, 
$\rho$ is the density, ${\bf v} = (u, v, w)$ is the velocity vector with Cartesian components $u$, $v$, and $w$, $p$ is the static pressure, 
and $H =   \gamma p / \{ \rho (\gamma-1) \} +  {\bf v} \cdot {\bf v}  / 2  $ is the specific total enthalpy with $\gamma=1.4$ (air), which is nondimensionalized by free 
stream values as described in Ref.~\cite{fun3d_manual:NASATM2016}, and closed by the nondimensionalized ideal gas law: 
\begin{eqnarray}
p = \rho T /  \gamma, 
\end{eqnarray}
where $T$ is the static temperature. Viscous terms are beyond the scope of the present work, but extensions to viscous-flow problems are possible through the hyperbolic diffusion system \cite{nishikawa_onetwothree_diffusion:JCP2014} and the hyperbolic Navier-Stokes systems \cite{liu_nishikawa_aiaa2016-3969,nishikawa_hyperbolic_ns:AIAA2015,NakashimaWatanabeNishikawa_AIAA2016-1101,LiLouNishikawaLuo_HNSrDG:JCP2021}, to which all the techniques proposed here for the Euler equations are directly applicable.

\begin{figure}[t] 
\begin{center}
\begin{minipage}[b]{0.75\textwidth}
\begin{center} 
          \includegraphics[width=0.4\textwidth,trim=0 0 0 0,clip]{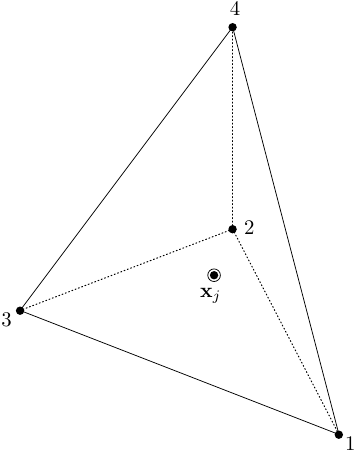}
        \caption{ A tetrahedral cell $j$ with a cell center at ${\bf x}_j$.} 
 \label{tetra_cell_ccfv}
\end{center}
\end{minipage}
\end{center}
\end{figure}

{\color{black}

On an unstructured tetrahedral grid, we store numerical solutions at cells and approximate the integral form of the conservation law over a cell $j$ (see Figure \ref{tetra_cell_ccfv} for an example of a tetrahedral cell): 
\begin{eqnarray}
\int_{V_j}  \frac{ \partial {\bf u} }{\partial t}  \, dV   +\int_{V_j}  \mbox{div}{\cal F} \, dV = \int_{V_j}  {\bf s} \, dV, 
\label{int_form}
\end{eqnarray}
or
\begin{eqnarray}
\frac{ d {\bf \overline{u}}_j }{d t}   +  \frac{1}{V_j} \oint_{\partial V_j} {\bf f} \, ds =  {\bf \overline{s}}_j  , 
\label{int_form2}
\end{eqnarray}
where $V_j$ denotes the volume of the cell $j$, $\partial V_j$ denotes the boundary of the cell $j$, ${\bf f}  = {\cal F} \hat{\bf n}$ is the flux projected in the unit outward normal vector $\hat{\bf n}$, and ${\bf \overline{u}}_j$ and ${\bf \overline{s}}_j$ are the cell averages of the solution and the source term, respectively. Following the QUICK scheme \cite{Leonard_QUICK_CMAME1979} and deconvolution methods \cite{FeliceDenaroMoela:NHT1993,Denaro:IJNMF1996,Denaro_CMAME:2015}, we store a point-valued solution ${\bf u}_j$ at a cell center ${\bf x}_j$, not a cell-average ${\bf \overline{u}}_j$. This simplifies the computation of an accurate solution at a face for a flux evaluation, i.e., we interpolate a point value at a face from point values at cells by a quadratic polynomial, instead of reconstructing a point value from cell averages. Therefore, there is no need to preserve a cell average and thus the quadratic interpolation does not involve any geometric moments. It also simplifies a gradient calculation method in a similar manner, as we will discuss later. To achieve third-order spatial accuracy, we need to discretize the surface integral and the cell averages of the time derivative and the source term by high-order quadrature formulas. To simplify the implementation in a second-order finite-volume code, we will construct quadrature formulas in the form of a singe-point quadrature formula with a correction. This process is especially important for the time derivative because it converts the time derivative of the cell average $ d {\bf \overline{u}}_j  / {d t} $ to that of the point value stored at cell, $ d {\bf {u}}_j / {d t} $, so that we can evolve our numerical solutions in time. As we will show later, by selecting the geometric centroid as the cell center ${\bf x}_j$: 
\begin{eqnarray}
{\bf x}_j = 
 \frac{1}{V_j} \int_{V_j}  {\bf x}    \, dV =
 \frac{1}{4} \sum_{i=1}^4 {\bf x}_i , 
 \label{centroid_definition_important}
\end{eqnarray}
of a tetrahedral grid cell, we find that the correction is merely a curvature term involving second derivatives of the solution and fluxes. These second derivatives can be easily computed in a one-dimensional grid as in Ref.~\cite{Nishikawa_3rdQUICK:2020}, but not in two- and three-dimensional unstructured grids, requiring extra computing time and additional memory. Our goal is to eliminate the need to compute and store the second derivatives. This is possible by utilizing an efficient projected-derivative formula that we will discuss in the next section.
}


\section{Efficient Projected-Derivative Formula via Nodal Values}
\label{Efficient_Projected_Derivative_Formula_via_Nodal_Values}
\indent 
\indent
Before discussing a third-order accurate discretization, we recall an efficient projected-derivative formula \cite{Frink:AIAAJ1998} that exploits function values available at nodes, and we derive a slightly generalized formula, which will be used as the foundation on which we will construct a third-order accurate discretization. This is the key idea in the construction of an economical third-order cell-centered finite-volume discretization, as we will show later.

\begin{figure}[t] 
\begin{center}
\begin{minipage}[b]{0.95\textwidth}
\begin{center}
          \includegraphics[width=0.45\textwidth,trim=0 0 0 0,clip]{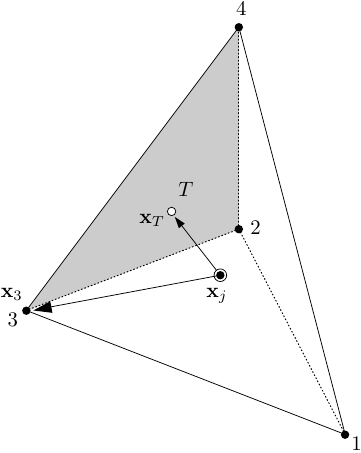}
        \caption{ Derivative projected to a face center ${\bf x}_T$ and a node $3$ from the cell centroid in a tetrahedral cell. Function values are assumed to be available 
        at the nodes.   }
 \label{efficient_projected_derivative}
\end{center}
\end{minipage}
\end{center}
\end{figure}

Consider an arbitrary function $g$ and a tetrahedral cell as shown in Figure \ref{efficient_projected_derivative}. If the function values are available at nodes $\{1,2,3,4\}$, the Green-Gauss gradient over the cell $j$ is given by 
\begin{eqnarray}
\nabla g_j^{GG} = \frac{1}{V_j} \sum_{i=1}^4  \overline{g}_i (-{\bf n}_i)  = \frac{1}{3 V_j} \sum_{i=1}^4 g_i {\bf n}_i  , 
\end{eqnarray}
where ${\bf n}_i$ is the inward normal vector of the triangular face opposite to a node $i$, whose magnitude is equal to the area of the triangular face, $ \overline{g}_i$ is the arithmetic average of the values of $g$ stored at the three nodes of that face, $V_j$ is the volume of the tetrahedral cell, $g_i$ is the value of $g$ at the node $i$. If ${\bf x}_j$ is the geometric centroid of the cell $j$ and ${\bf x}_{T}$ is the geometric centroid of a triangular face $T$ of the cell: e.g., for the triangular face $T$ opposite to the node $1$ in Figure \ref{efficient_projected_derivative} (shaded),
\begin{eqnarray} 
{\bf x}_j = \frac{ {\bf x}_1 + {\bf x}_2 + {\bf x}_3 + {\bf x}_4 }{4}
, \quad
{\bf x}_T = \frac{  {\bf x}_2 + {\bf x}_3 + {\bf x}_4 }{3}, 
\end{eqnarray}
we have
\begin{eqnarray}
{\bf n}_1 \cdot ( {\bf x}_{T} - {\bf x}_j ) = - \frac{3 V_j}{4}, 
\quad
{\bf n}_2 \cdot ( {\bf x}_{T} - {\bf x}_j ) = {\bf n}_3 \cdot ( {\bf x}_{T} - {\bf x}_j ) = {\bf n}_4 \cdot ( {\bf x}_{T} - {\bf x}_j ) = \frac{3 V_j}{12}.
\end{eqnarray}
Then, the Green-Gauss gradient projected along the line segment from ${\bf x}_j$ to ${\bf x}_{T}$ can be simplified as  
\begin{eqnarray}
  \nabla g_j^{GG}  \cdot ( {\bf x}_{T} - {\bf x}_j )
 =
   \frac{1}{4} \left(  \frac{g_2 + g_3 + g_4 }{3} - g_1 \right),
 \label{simplified_nablag_dx}
\end{eqnarray}
which results in the following linear reconstruction formula from the cell centroid to the face centroid: 
\begin{eqnarray}
 g_j + \nabla g_j^{GG} \cdot ( {\bf x}_{T} - {\bf x}_j )
 =
 g_j +  \frac{1}{4} \left(  \frac{g_2 + g_3 + g_4 }{3} - g_1 \right). 
 \label{simplified_nablag_dx_usm3d}
\end{eqnarray}
This is an efficient formula that does not require explicit computations and storage for the solution gradient; it was proposed originally for, and is still used in, the USM3D code \cite{Frink:AIAAJ1998,pandya_etal:AIAAJ2016}.

The projected-derivative formula (\ref{simplified_nablag_dx}) can be generalized to an arbitrary point ${\bf x}$. Notice that the formula (\ref{simplified_nablag_dx}) can be viewed as a finite-difference approximation between ${\bf x}_j$ and ${\bf x}_{T}$: 
\begin{eqnarray}
\nabla g_j ^{GG} \cdot ( {\bf x}_{T} - {\bf x}_j ) 
=
 \frac{1}{4} \left(  \frac{g_2 + g_3 + g_4 }{3} - g_1 \right)
= \left( \frac{ \overline{g}_{T} - \overline{g}_j  }{  |   {\bf x}_{T} - {\bf x}_j    |  }  \right)  |   {\bf x}_{T} - {\bf x}_j    |   = \overline{g}_{T} - \overline{g}_j ,
\label{generalized_LRviaNV00}
\end{eqnarray} 
where $\overline{g}_{T}$ and $\overline{g}_{j}$ are the function values linearly interpolated at the geometric centroids of the triangular face and the tetrahedron, respectively: 
\begin{eqnarray}
\overline{g}_{T}  =   \frac{ g_2 + g_3 + g_4 }{3},
\quad
\overline{g}_j   =   \frac{ g_1 + g_2 + g_3 + g_4 }{4}.
\label{gjk_gj_geometric_centroids_of_tetra}
\end{eqnarray}
Then, we have the following generalization applicable to an arbitrary point ${\bf x}$:
\begin{eqnarray}
\nabla g^{GG} \cdot ( {\bf x} - {\bf x}_j )   =   {g} - \overline{g}_j, 
\label{generalized_LRviaNV00_general}
\end{eqnarray} 
where $g$ is a linearly-exact value at ${\bf x}$ (e.g., a linearly interpolated value). The generalized formula, Eq. (\ref{generalized_LRviaNV00_general}), can be used not only at the face center as in Eq. (\ref{generalized_LRviaNV00}) with $g=\overline{g}_{T}$ but also at a node $i$ with $g={g}_i$, 
\begin{eqnarray}
\nabla g_j^{GG} \cdot ( {\bf x}_{i} - {\bf x}_j )  = {g}_{i} - \overline{g}_j ,
 \label{generalized_LRviaNV00_m}
\end{eqnarray}
where $g_{i}$ is a linearly-exact value of $g$ at the node $i$. 
For example, as indicated in Figure \ref{efficient_projected_derivative}, we have, at the node $3$, 
\begin{eqnarray}
\nabla g_j^{GG} \cdot ( {\bf x}_{3} - {\bf x}_j )  =   {g}_{3} - \overline{g}_j .
\end{eqnarray}
{\color{black} As one would expect, an equivalent formula exists for a triangle: 
\begin{eqnarray}
\nabla g_T^{GG} \cdot ( {\bf x}_{i} - {\bf x}_T )  =   {g}_{i} - \overline{g}_T, 
 \label{generalized_LRviaNV00_m_2D}
\end{eqnarray}
where $\nabla g_T^{GG}$ denotes the Green-Gauss gradient defined over a triangle $T$ and $i$ is one of its vertices (e.g., $i = 2, 3 $, or $4$ for the triangular face $T$ in Figure \ref{efficient_projected_derivative}).} Notice that these formulas allow us to express the gradient of a function in terms of nodal values only. In this paper, we will employ these formulas to express a quadratic term involving a Hessian matrix in terms of the nodal gradients, thereby eliminating the need to explicitly compute and store second derivatives. This is how we eliminate second derivatives from a third-order accurate cell-centered discretization as we discuss in detail in the next section. 

\section{Economical Third-Order Cell-Centered Discretization}
\label{Economical_Third_Order_Cell_Centered_Discretization}
 
In this section, we describe a new economical third-order cell-centered discretization constructed with the efficient projected-derivative formula. There are three key items. The first is a quadratically exact flux quadrature formula with a single numerical flux per face, where we only need one numerical flux per face but need to add a curvature correction term for quadratic exactness. Yet, we can eliminate all second derivatives from the curvature term by utilizing the efficient projected-derivative formula derived in the previous section. {\color{black} The second is a quadratically exact source quadrature formula expressed in the form of a point value at the cell centroid with a curvature correction term, where again second derivatives are eliminated by the efficient projected-derivative formula. The third} is that the gradient computation at nodes from solutions at cells, which requires a quadratic fit, can be formulated in a simple manner because the solutions at cells and the gradients at nodes are both point values. In addition, since there is no need to compute second derivatives, the resulting LSQ gradient formula can be written in the same form as a linear LSQ gradient formula and thus the cost of gradient computations will be the same if the number of neighbors are the same in linear and quadratic methods. These key items allow an efficient quadratic interpolation scheme to be constructed that does not require second derivatives at all. Then, we will discuss the implementation of weak boundary conditions. Finally, we will discuss how the resulting discretization is more efficient than a conventional third-order discretization.

\subsection{{\color{black} Flux Quadrature}}
\indent 
\indent 
Consider a tetrahedral cell $j$ and its face neighbors $\{ k_j \}$. Let $T_{jk}$ be the triangular face shared with a neighbor cell $k \in \{ k_j \}$. Then, the conservative integral form, Eq. (\ref{int_form2}), is written as 
\begin{eqnarray}
\frac{ d  {\color{black} \overline{\bf u}_j } }{d t}   +  \frac{1}{V_j}      \sum_{ k \in \{  k_j \}  }   \int_{T_{jk} }   {\bf f} \,  ds = {\color{black} \overline{\bf s}_j } , 
\end{eqnarray}
where $ {\bf f}  =  {\cal F}   \hat{\bf n}_T$ and $\hat{\bf n}_T$ is the unit outward normal vector of the face $T_{jk}$. To achieve third-order accuracy, we need to discretize the surface integral by a quadratically-exact quadrature method. An efficient quadrature formula can be derived by integrating a quadratic flux over the face: 
\begin{eqnarray}
 \int_{T_{jk} }   {\bf f} \,  ds 
 =
 \int_{T_{jk} }     \left[    {\bf f}( {\bf x}_T) +   \nabla  {\bf f}_T  \cdot ( {\bf x} - {\bf x}_{T} )  + \frac{1}{2}  ( {\bf x} - {\bf x}_{T} )^t 
    \nabla^2 {\bf f}_T 
   ( {\bf x} - {\bf x}_{T} )          \right]  ds ,
   \label{surface_quadrature_00}
\end{eqnarray}
where the superscript $t$ indicates transpose and $ {\bf x}_{T}$ is the geometric centroid of the triangular face. The linear term vanishes on integration over the triangle and we are left with
 \begin{eqnarray}
 \int_{ T_{jk} }{\bf f} \,  ds 
 =
  {\bf f}( {\bf x}_T)   | {\bf n}_T |
  +
  \frac{1}{2}    \int_{ T_{jk} }  ( {\bf x} - {\bf x}_{T} )^t 
    \nabla^2 {\bf f}_T 
   ( {\bf x} - {\bf x}_{T} )      ds ,
   \label{surface_quadrature_01}
\end{eqnarray}
where ${\bf n}_T$ is the outward normal vector of the face and $ | {\bf n}_T |$ is the area of the face. The quadratic terms can be analytically integrated over a triangle using local coordinates as described in Refs.~\cite{EfficientIntegrationTriaTetra:2017,LiuVinokur:JCP1997}:
 \begin{eqnarray}
  \int_{ T_{jk} }   (x - x_T )^2    ds   =  \frac{1}{12}    \sum_{i=1}^3    (x_i - x_T )^2    | {\bf n}_T |,
  &   &
    \int_{ T_{jk} }   (x - x_T )(y - y_T )    ds  =  \frac{1}{12}    \sum_{i=1}^3    (x_i - x_T )(y_i - y_T )    | {\bf n}_T |,    \\ [2ex]
    \int_{ T_{jk} }   (y - y_T )^2    ds  = \frac{1}{12}    \sum_{i=1}^3    (y_i - y_T )^2    | {\bf n}_T |,  
      &  &
        \int_{ T_{jk} }   (y - y_T )(z - z_T )    ds    =   \frac{1}{12}    \sum_{i=1}^3    (y_i - y_T )(z_i - z_T )    | {\bf n}_T |,   \\ [2ex] 
    \int_{ T_{jk} }   (z - z_T )^2    ds   =   \frac{1}{12}    \sum_{i=1}^3    (z_i - z_T )^2    | {\bf n}_T |,  
  & & 
    \int_{ T_{jk} }   (z - z_T )(x - x_T )   ds    =   \frac{1}{12}    \sum_{i=1}^3    (z_i - z_T )  (x_i - x_T )   | {\bf n}_T |, 
\end{eqnarray} 
and therefore, Eq. (\ref{surface_quadrature_01}) can be written as
 \begin{eqnarray}
 \int_{T_{jk} } {\bf f} \,  ds 
 =
  {\bf f}( {\bf x}_T)   | {\bf n}_T |
  +
  \frac{1}{24}  \sum_{i=1}^3    ( {\bf x}_i - {\bf x}_{T} )^t 
    \nabla^2 {\bf f}_T 
   ( {\bf x}_i - {\bf x}_{T} )     | {\bf n}_T |, 
\end{eqnarray}
where $i= 1$, $2$, and $3$ are the nodes of the triangular face $T_{jk}$. It can be further simplified by the efficient projected-gradient formula, Eq. {\color{black} (\ref{generalized_LRviaNV00_m_2D})}. Noting that Eq. {\color{black} (\ref{generalized_LRviaNV00_m_2D})} can be expressed in terms of the directional $\nabla$: $\left[ ( {\bf x}_{i} - {\bf x}_T ) \cdot  \nabla^{GG} \right] g_T = {g}_{i} - \overline{g}_T $, we find
 \begin{eqnarray}
 \int_{T_{jk} } {\bf f} \,  ds 
 &=&
  {\bf f}( {\bf x}_T)   | {\bf n}_T |
  +
  \frac{1}{24}  \sum_{i=1}^3    ( {\bf x}_i - {\bf x}_{T} )^t 
    \nabla^2 {\bf f}_T 
   ( {\bf x}_i - {\bf x}_{T} )     | {\bf n}_T |   
  \nonumber \\ [2ex]  
  &=&
  {\bf f}( {\bf x}_T)   | {\bf n}_T |
  +
  \frac{1}{24}  \sum_{i=1}^3   \left[ 
   ( {\bf x}_i - {\bf x}_{T} ) \cdot 
\nabla^{GG} \right]  \nabla {\bf f}_T
   ( {\bf x}_i - {\bf x}_{T} )     | {\bf n}_T |   
  \nonumber \\ [2ex] 
  &=&
  {\bf f}( {\bf x}_T)   | {\bf n}_T |
  +
  \frac{1}{24}  \sum_{i=1}^3     (     \nabla {\bf f}_i -   \nabla {\bf f}_T  )   ( {\bf x}_i - {\bf x}_{T} )  | {\bf n}_T |   
  \nonumber \\ [2ex] 
& =&
   \left[  
  {\bf f}( {\bf x}_T)   
  +
  \frac{1}{24}     \sum_{i=1}^3  (     \nabla {\bf f}_i -   \nabla {\bf f}_T  )  ( {\bf x}_i - {\bf x}_{T} )     \right]   | {\bf n}_T |  ,
  \label{flux_corr_simplify_00}
\end{eqnarray}
and furthermore,
 \begin{eqnarray}
 \int_{T_{jk} } {\bf f} \,  ds 
& =&
   \left[  
  {\bf f}( {\bf x}_T)   
  +
  \frac{1}{24}     \sum_{i=1}^3  (     \nabla {\bf f}_i -   \nabla {\bf f}_T  )  ( {\bf x}_i - {\bf x}_{T} )     \right]   | {\bf n}_T | 
  \nonumber \\ [2ex]
& =&
   \left[  
  {\bf f}( {\bf x}_T)   
  +
  \frac{1}{24}     \sum_{i=1}^3        \nabla {\bf f}_i     ( {\bf x}_i - {\bf x}_{T} )    
-
  \frac{1}{24}     \sum_{i=1}^3    \nabla {\bf f}_T  ( {\bf x}_i - {\bf x}_{T} )     \right]   | {\bf n}_T | 
  \nonumber \\ [2ex]
& =&
   \left[  
  {\bf f}( {\bf x}_T)   
  +
  \frac{1}{24}     \sum_{i=1}^3        \nabla {\bf f}_i     ( {\bf x}_i - {\bf x}_{T} )   
  -
  \frac{1}{24}    \nabla {\bf f}_T   \sum_{i=1}^3    
  ( {\bf x}_i - {\bf x}_{T} )  
       \right]   | {\bf n}_T | 
  \nonumber \\ [2ex]
& =&
   \left[  
  {\bf f}( {\bf x}_T)   
  +
  \frac{1}{24}     \sum_{i=1}^3        \nabla {\bf f}_i     ( {\bf x}_i - {\bf x}_{T} )   
  -
  \frac{1}{24}   \nabla {\bf f}_T   \left( 
   \sum_{i=1}^3    
   {\bf x}_i - 3 {\bf x}_{T}  
  \right)
       \right]   | {\bf n}_T | 
  \nonumber \\ [2ex]
& =&
   \left[  
  {\bf f}( {\bf x}_T)   
  +
  \frac{1}{24}     \sum_{i=1}^3        \nabla {\bf f}_i     ( {\bf x}_i - {\bf x}_{T} )   
  -
  \frac{1}{24}   \nabla {\bf f}_T   \left( 
 3 {\bf x}_{T}   - 3 {\bf x}_{T}  
  \right)
       \right]   | {\bf n}_T | 
  \nonumber \\ [2ex]
& =&
 \left[  
  {\bf f}( {\bf x}_T)   
  +
  \frac{1}{24}     \sum_{i=1}^3     \nabla {\bf f}_i  ( {\bf x}_i - {\bf x}_{T} )    \right]   | {\bf n}_T |.
  \label{flux_corr_simplify_01}
\end{eqnarray}
Then, the flux at the face centroid $ {\bf f}( {\bf x}_T) $ is replaced by a numerical flux $ \Phi_{jk}$:
 \begin{eqnarray}
 \int_{T_{jk} } {\bf f} \,  ds 
=
 \left[  
 \Phi_{jk}
  +
  \frac{1}{24}     \sum_{i=1}^3     \nabla {\bf f}_i   ( {\bf x}_i - {\bf x}_{T} )    \right]   | {\bf n}_T |.
\end{eqnarray}
For the numerical flux $ \Phi_{jk}$, we employ the Roe flux \cite{Roe_JCP_1981} in this work, unless otherwise stated: 
\begin{eqnarray}
 \Phi_{jk}
=   \frac{1}{2} \left[  {\bf f}( {\bf w}_L) + {\bf f}( {\bf w}_R)    \right]
-
 \frac{1}{2}  {\bf D} 
\left(   {\bf w}_R -  {\bf w}_L \right)
, 
\label{numerical_flux_Phijk}
\end{eqnarray} 
where $ {\bf D}$ is a dissipation matrix evaluated at the Roe-averages, and ${\bf w}_L$ and ${\bf w}_R$ are the primitive variables ${\bf w} = (\rho, {\bf v}, p)$ quadratically interpolated at the face centroid from the cell $j$ and the neighbor cell $k$, respectively. Finally, we avoid directly computing the flux gradients, which can be expensive, by using the chain rule:
 \begin{eqnarray}
 \int_{T_{jk} } {\bf f} \,  ds 
=
 \left[  
 \Phi_{jk}
  +
  \frac{1}{24}     \sum_{i=1}^3   \left(  \frac{\partial {\bf f}}{\partial {\bf w}} \right)_{\!\! i}  \!\!    \nabla {\bf w}_i   ( {\bf x}_i - {\bf x}_{T} ) \right]  | {\bf n}_T |,
  \label{flux_quadrature_CQ}
\end{eqnarray}
where $ \nabla {\bf w}_i $ is a solution gradient computed and stored at a node $i$ and the normal flux Jacobian is evaluated with a solution reconstructed at the node $i$ {\color{black} but with the same face normal $\hat{\bf n}_T$ (i.e., $ {\bf f}  =  {\cal F} \hat{\bf n}_T$) for any $i$}. {\color{black} Note that the resulting quadrature formula is in the form of the numerical flux at the face centroid with a correction term, as desired, and the correction term does not involve second derivatives. In the next section, we derive a source term quadrature formula, which will be applied to the time derivative term as well.}

{\color{black} 
\subsection{{\color{black} Source Term Quadrature}}
\indent 
\indent 
Consider the cell-averaged source term: 
 \begin{eqnarray}
 \overline{\bf s}_j 
 =
 \frac{1}{V_j} 
 \int_{V_j}  {\bf s} \, dV.
\end{eqnarray}
We wish to derive a quadratically-exact quadrature formula in the form of a point value ${\bf s}_j = {\bf s} ({\bf x}_j)$ at the cell center ${\bf x}_j$ with a correction. Expanding the source term up to the quadratic term around the centroid ${\bf x}_j$:
\begin{eqnarray}
 \overline{\bf s}_j 
 =
 \frac{1}{V_j} 
 \int_{V_j} 
   \left[    {\bf s}_j   +   \nabla  {\bf s}_j  \cdot ( {\bf x} - {\bf x}_{j} )  + \frac{1}{2}  ( {\bf x} - {\bf x}_{j} )^t 
    \nabla^2 {\bf s}_j 
   ( {\bf x} - {\bf x}_{j} )          \right]  
  \, dV,
\end{eqnarray}
we obtain
\begin{eqnarray}
 \overline{\bf s}_j 
 =
  {\bf s}_j  
  +  
 \frac{1}{2 V_j} 
 \int_{V_j} 
 ( {\bf x} - {\bf x}_{j} )^t 
    \nabla^2 {\bf s}_j 
   ( {\bf x} - {\bf x}_{j} )   
  \, dV.
\end{eqnarray}
The quadratic terms can be analytically integrated over a tetrahedron \cite{EfficientIntegrationTriaTetra:2017,LiuVinokur:JCP1997} as
 \begin{eqnarray}
  \int_{ V_j }   (x - x_j )^2    dV   =  \frac{1}{20}    \sum_{i=1}^4    (x_i - x_j )^2    V_j ,
  &   &
    \int_{V_j }   (x - x_j )(y - y_j )    dV  =  \frac{1}{20}    \sum_{i=1}^4    (x_i - x_j )(y_i - y_j )    V_j ,   \\ [2ex]
    \int_{ V_j }   (y - y_j )^2    dV  = \frac{1}{20}    \sum_{i=1}^4    (y_i - y_j )^2   V_j ,
      &  &
        \int_{ V_j}   (y - y_j )(z - z_j )    dV    =   \frac{1}{20}    \sum_{i=1}^4   (y_i - y_j )(z_i - z_j ) V_j ,  \\ [2ex] 
    \int_{ V_j }   (z - z_j )^2    dV   =   \frac{1}{20}    \sum_{i=1}^4    (z_i - z_j )^2    V_j ,
  & & 
    \int_{V_j}   (z - z_j )(x - x_j )   dV    =   \frac{1}{20}    \sum_{i=1}^4   (z_i - z_j )  (x_i - x_j ) V_j ,
\end{eqnarray} 
where i = 1, 2, 3, and 4 are the nodes of the tetrahedron $V_j$, and therefore, 
\begin{eqnarray}
 \overline{\bf s}_j 
 =
  {\bf s}_j  
  + 
   \frac{1}{40}
   \sum_{i=1}^4
    ( {\bf x}_i - {\bf x}_{j} )^t 
    \nabla^2 {\bf s}_j 
   ( {\bf x}_i - {\bf x}_{j} ),
\end{eqnarray}
which can be simplified, in a similar manner as described in Eqs.~(\ref{flux_corr_simplify_00}) and (\ref{flux_corr_simplify_01}), by the efficient projected-gradient formula, Eq. {\color{black} (\ref{generalized_LRviaNV00_m})}, as
\begin{eqnarray}
 \overline{\bf s}_j 
& = &
  {\bf s}_j  
  + 
   \frac{1}{40}
   \sum_{i=1}^4
  (  \nabla {\bf s}_i -  \nabla {\bf s}_j 
  )
   ( {\bf x}_i - {\bf x}_{j} )
 \nonumber \\ [2ex]
 &=&
  {\bf s}_j 
  + 
   \frac{1}{40}
   \sum_{i=1}^4
\nabla {\bf s}_i  ( {\bf x}_i - {\bf x}_{j} ).
\end{eqnarray}
Therefore, the source term quadrature formula is expressed, as desired, in the form of a point value at the centroid ${\bf s}_j$ with a correction, and the second derivatives have been eliminated from the correction term. In the same manner, this source quadrature formula is directly applicable to the time derivative term for unsteady problems: 
\begin{eqnarray} 
  \frac{ \partial \overline{\bf u}_j }{\partial t} 
=
 \frac{ \partial {\bf u}_j }{\partial t}  
  + 
   \frac{1}{40}
   \sum_{i=1}^4
\nabla 
\left(  \frac{ \partial {\bf u}_i }{\partial t}  \right)  ( {\bf x}_i - {\bf x}_{j} ).
\end{eqnarray}
The final form of the third-order discretization is given by
\begin{eqnarray}  
\frac{ \partial {\bf u}_j }{\partial t}  
  +
  \frac{1}{V_j}
     \sum_{ k \in \{  k_j \}  }  
   \left[  
 \Phi_{jk}
  +
  \frac{1}{24}     \sum_{i=1}^3   \left(  \frac{\partial {\bf f}}{\partial {\bf w}} \right)_{\!\! i}  \!\!    \nabla {\bf w}_i   ( {\bf x}_i - {\bf x}_{T} ) \right]  | {\bf n}_T |
  =  
  {\bf s}_j 
  + 
   \frac{1}{40}
   \sum_{i=1}^4
\nabla 
\left[ {\bf s}_i 
 - \left(  \frac{ \partial {\bf u}_i }{\partial t}  \right) 
\right] ( {\bf x}_i - {\bf x}_{j} ).
\label{final_third_order_res}
\end{eqnarray}
Each correction term is a curvature term not involving second derivatives, where the gradients need to be available at nodes. This is aligned well with the main idea of the FANG approach. Inspired by the work of Zhang \cite{zhang_etal:2015,Zhang:Preprint2017,zhang:phd}, where he proposed to compute gradients at nodes from solution values at cells and obtain the cell gradient needed in a cell-centered finite-volume method by averaging nodal gradients over cell nodes, we developed the FANG approach \cite{NishikawaWhite_FANG:jcp2020}, where solution interpolation is performed with a common gradient defined at a face by averaging nodal gradients over the face (rather than over the cell). It brings substantial saving in computing time and memory for gradients on triangular and tetrahedral grids, leads to smaller residual stencils, and allows efficient parallel implementations \cite{NishikawaWhite_FANG:jcp2020,Nishikawa_FANG_AQ:Aviation2020,WhiteNishikawaBaurle_scitech2020,NishikawaWhite:AIAA2021-2720}. In this work, we extend the FANG approach to third-order accuracy based on the above discretization. For third-order accuracy, the nodal gradient method must be exact for quadratic functions, and also the solution interpolation needs to be upgraded to a quadratic interpolation. Normally, a quadratic interpolation would involve second derivatives of the primitive variables, but we can eliminate them all together by utilizing the efficient projected-derivative formula just like we have done in deriving the third-order discretization in the above. The resulting interpolation scheme will be employed not only for computing ${\bf w}_L$ and ${\bf w}_R$ at the face centroid but also for interpolating the solution at nodes to evaluate the flux Jacobian in Eq. (\ref{flux_quadrature_CQ}). 
}
 
\subsection{Computing Gradients at Nodes by Quadratic Least-Squares Fit}
\label{LSQ_stencils}
\indent
\indent
Before deriving an efficient {\color{black} interpolation} scheme, we briefly describe the method for computing gradients at nodes from the solution values stored at cells. 
Consider a quadratic polynomial for a scalar variable $w$ centered at a node $i$: 
\begin{eqnarray}
w({\bf x}) = w_i +  {\nabla}  w_i \cdot ( {\bf x}- {\bf x}_i ) + \frac{1}{2}  ( {\bf x} - {\bf x}_i)^t  {\nabla}^2  w_i  ( {\bf x} - {\bf x}_i ),
\end{eqnarray}
where $w_i$, $ {\nabla} w_i$, and $ {\nabla}^2 w_i $ are unknown solution value, gradient, and Hessian matrix at the node $i$, and ${\bf x}_i$ is the position of the node $i$. Note that the numerical solutions at cells are point values and there is no need to impose the condition that a cell average of the polynomial is equal to $w_i$ (as typically required in other finite-volume methods). Then, we fit the above polynomial directly over the values stored at a set $\{ k_i \}$ of neighbor cells and obtain an $N_{\{ k_i \}}^{LSQ}$$\times$$10$ system, where $N_{\{ k_i \}}^{LSQ}$ is the number of neighbors in the set $\{ k_i \}$: the $k$-th row of the resulting system is given by 
\begin{eqnarray} 
 w_i +  {\nabla}  w_i \cdot ( {\bf x}_k - {\bf x}_i ) + \frac{1}{2}  ( {\bf x}_k  - {\bf x}_i )^t  {\nabla}^2  w_i  ( {\bf x}_k - {\bf x}_i )   =   w_k,
 \quad
 k \in  \{ k_i \}, 
\end{eqnarray}
where ${\bf x}_k$ is the centroid of the $k$-th neighbor cell and $w_k$ is the corresponding solution value stored. It is possible to construct a weighted LSQ problem, e.g., using the inverse-distance weighting, but we do not apply any weighting in this work.  Assuming $N_{\{ k_i \}}^{LSQ}  \ge 10$, which is required for a well posed quadratic fit, we solve the system by the QR factorization via the Householder transformation \cite{GilbertStrang_1980}, and obtain 
\begin{eqnarray}
{\nabla}  w_i
=  
\sum_{k \in \{ k_i \} } 
\left[
\begin{array}{c}
c_{ik}^x  \\ [1.2ex]
c_{ik}^y \\ [1.2ex]
c_{ik}^z
\end{array}
\right] 
w_k  , 
\end{eqnarray}
where $c_{ik}^x$, $c_{ik}^y$, $c_{ik}^z$ are the LSQ coefficients to be computed and stored at all nodes once for a given stationary grid. Note that we also obtain the coefficients for the solution interpolation for $w_i$ as well as those for the second derivatives at nodes but we do not store them since, as we will discuss later, the second derivatives will be eliminated from the quadratic {\color{black} interpolation} scheme. For the solution value, $w_i$, which is required for the flux Jacobian evaluation, we would rather compute it by using the quadratic {\color{black} interpolation} scheme instead of the LSQ method, which would make it easier to handle unrealizable states. 

For the choice of the neighbors $\{ k_i \}$, as discussed in Refs.~\cite{NishikawaWhite:AIAA2021-2720,VULCAN_FANGplus_scitech2022}, it is known that adding extra neighbors to the basic set $\{   k_i^{(0)}  \}$ of cells sharing a node $i$ helps make the solver more robust. For example, we may add those sharing at least one face with any of $\{   k_i^{(0)}  \}$ (face neighbors) or those sharing at least one node with any of $\{   k_i^{(0)}  \}$ (node neighbors) to $\{   k_i^{(0)}  \}$, which respectively lead to the methods referred to as FANG+fn and FANG+nn \cite{NishikawaWhite:AIAA2021-2720,VULCAN_FANGplus_scitech2022}. In this work, we employ FANG+nn unless otherwise stated, and therefore, $\{ k_i \}$ is the sum of the cells sharing a node and their node neighbors. In all test cases considered in this work, we have $N_{\{ k_i \}}^{LSQ}  \ge 10$ and the quadratic LSQ problem has been found to be well posed.

\subsection{Efficient Quadratic {\color{black} Interpolation} via Nodal Gradients}
 
\begin{figure}[t] 
\begin{center}
\begin{minipage}[b]{0.8\textwidth}
\begin{center}
          \includegraphics[width=0.65\textwidth,trim=0 0 0 0,clip]{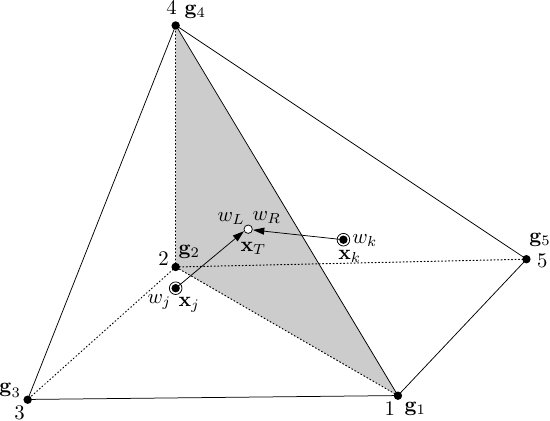}
        \caption{A face between two tetrahedral cells, $j$ and $k$: $w_j$ and $w_k$ are the solutions stored at the centroids ${\bf x}_j$ and ${\bf x}_k$; ${\bf g}_1$, ${\bf g}_2$, ${\bf g}_3$, ${\bf g}_4$, and ${\bf g}_5$ are the gradients stored at nodes.         
        As an example, the solutions interpolated at the face centroid from the left and right cells are indicated as $w_L$ and $w_R$, respectively.   }
       \label{tetra_quad_reconstruction_j_k}
\end{center}
\end{minipage}
\end{center}
\end{figure}

We now derive an efficient quadratic {\color{black} interpolation} scheme that does not require second derivatives. 
 For simplicity, let us consider a scalar variable $w$ and a quadratic polynomial over a cell $j$:
\begin{eqnarray}
w({\bf x}) = w_j +  {\nabla}  w_j \cdot ( {\bf x}- {\bf x}_j ) + \frac{1}{2}  ( {\bf x} - {\bf x}_j )^t   {\nabla}^2 w_j ( {\bf x} - {\bf x}_j ), 
\label{quadratic_reconstruction_00}
\end{eqnarray}
where $w_j$ is a solution value stored at the cell $j$, $ {\nabla} w_j$ is the gradient at a cell center ${\bf x}_j$, and $ {\nabla}^2 w_j $ is the Hessian matrix of $w$. It is emphasized that the solution value $w_j$ stored at the cell is a point value solution at ${\bf x}_j$, not a cell average. Then, the left state at an arbitrary point, ${\bf x}$, between the cells $j$ and $k$ (e.g., at a node or at the face centroid) is evaluated as
\begin{eqnarray}
w_L = w_j + {\nabla} w_j \cdot ( {\bf x} - {\bf x}_j ) + \frac{1}{2}  ( {\bf x} - {\bf x}_j )^t  {\nabla}^2 w_j  ( {\bf x} - {\bf x}_j ). 
\label{quadratic_reconstruction_01}
\end{eqnarray} 
See Figure \ref{tetra_quad_reconstruction_j_k}, which illustrates a face between $j$ and $k$, and the face centroid point ${\bf x}_T$. In the FANG approach, we employ $ {\nabla} w$ defined at a face, instead of ${\nabla} w_j$, in the solution {\color{black} interpolation}; this preserves second-order accuracy \cite{NishikawaWhite_FANG:jcp2020}. For third-order accuracy, we need to linearly extrapolate the gradient from the cell center to the face: 
\begin{eqnarray}
 \nabla w 
 =
    \nabla w_j  + {\nabla}^2 w_j  (  {\bf x}    -    {\bf x}_j  ). 
\label{linear_extrapolation_nablaw_j}
\end{eqnarray}
Substituting Eq. (\ref{linear_extrapolation_nablaw_j}) into Eq. (\ref{quadratic_reconstruction_01}), we obtain the following quadratic {\color{black} interpolation} formula:
\begin{eqnarray}
w_L 
&=&
w_j + 
\nabla w \cdot ( {\bf x}  - {\bf x}_j )   -   ( {\bf x}  - {\bf x}_j )^t  {\nabla}^2 w_j   (  {\bf x}    -    {\bf x}_j  ) +   \frac{1}{2}  ( {\bf x}  - {\bf x}_j )^t   {\nabla}^2 w_j    (   {\bf x}  - {\bf x}_j )  
\nonumber \\ [1ex]
&=&
w_j + \nabla w  \cdot ( {\bf x}  - {\bf x}_j ) - \frac{1}{2}  ( {\bf x}  - {\bf x}_j )^t   {\nabla}^2 w_j     (   {\bf x}  - {\bf x}_j )  
, 
\label{quadratic_reconstruction_02}
\end{eqnarray}
where $ \nabla w $ is computed by averaging the gradients stored at nodes of the face (i.e., linear interpolation), 
and similarly for the right state:
\begin{eqnarray}
w_R 
=
w_k + \nabla w  \cdot ( {\bf x}  - {\bf x}_k ) - \frac{1}{2}  ( {\bf x}  - {\bf x}_k )^t  {\nabla}^2 w_k     (   {\bf x} - {\bf x}_k )  .
\label{quadratic_reconstruction_02_wR}
\end{eqnarray}
Note that the location where the Hessian matrices $ {\nabla}^2 w_j$ and $ {\nabla}^2 w_k$ are evaluated is not important for quadratic {\color{black} interpolation} because the Hessian is constant for a quadratic function. However, we will use the values defined in the corresponding cells because it allows a dramatic simplification as we will show next.
 
In the FANG approach, the gradients are stored at nodes and we can simplify the curvature term of the quadratic {\color{black} interpolation} by applying the projected-derivative formula (\ref{generalized_LRviaNV00_general}). Let ${\bf g} = (g_x, g_y, g_z) = ( \partial_x w, \partial_y w, \partial_z w )$. Then, we find, by the formula (\ref{generalized_LRviaNV00_general}),
\begin{eqnarray}
{\nabla}^2 w_j   ( {\bf x} - {\bf x}_j ) 
= 
 \left[\begin{array}{c}
( \nabla {g}_x )_j^{GG} \cdot  ( {\bf x} - {\bf x}_j )   \\  [1.7ex] 
( \nabla {g}_y)_j^{GG} \cdot  ( {\bf x} - {\bf x}_j )   \\  [1.7ex] 
( \nabla {g}_z)_j^{GG} \cdot  ( {\bf x} - {\bf x}_j )  
\end{array}\right]
=
 \left[\begin{array}{c}
\displaystyle     {g}_x -  (\overline{g}_x)_j  \\  [1.7ex] 
\displaystyle     {g}_y -  (\overline{g}_y)_j  \\  [1.7ex] 
\displaystyle     {g}_z -  (\overline{g}_z)_j 
\end{array}\right] 
=
 {\bf g} -  \overline{\bf g}_j   
 ,
\end{eqnarray}
where $ \overline{\bf g}_j $ is the value of ${\bf g}$ linearly interpolated at ${\bf x}_j$, 
\begin{eqnarray}
\overline{\bf g}_j    = \frac{ {\bf g}_1 + {\bf g}_2 + {\bf g}_3 + {\bf g}_4 }{4} , 
\end{eqnarray}
and ${\bf g}$ of ${\bf g} -  \overline{\bf g}_j$ is the value linearly interpolated at ${\bf x}$.  
 Note that the projected-derivative formula (\ref{generalized_LRviaNV00_general}) is exact for linear functions and that is precisely what we need here for the second derivative terms because linear exactness for the gradient components is equivalent to quadratic exactness for the solution. Hence, the quadratic term is simplified as
\begin{eqnarray}
 \frac{1}{2}  ( {\bf x} - {\bf x}_j )^t {\nabla}^2 w_j   ( {\bf x} - {\bf x}_j )
=
 \frac{1}{2}  
( {\bf g} -  \overline{\bf g}_j )    \cdot   ( {\bf x} - {\bf x}_j ).
\end{eqnarray}
Therefore, the quadratic {\color{black} interpolation} scheme, Eq. (\ref{quadratic_reconstruction_02}), can be simplified as
\begin{eqnarray}
w_L  
&=& 
w_j + \nabla w    \cdot ( {\bf x} - {\bf x}_j ) - \frac{1}{2}  
( {\bf g} -  \overline{\bf g}_j )    \cdot   ( {\bf x}  - {\bf x}_j ) 
\nonumber \\ [1ex]
&=& 
w_j +  {\bf g}  \cdot ( {\bf x} - {\bf x}_j ) -
 \frac{1}{2} ( {\bf g}  -  \overline{\bf g}_j )    \cdot   ( {\bf x} - {\bf x}_j )
 \nonumber \\ [1ex]
&=& 
w_j +  \frac{1}{2}    \left(   {\bf g} +  \overline{\bf g}_j   \right)  \cdot   ( {\bf x} - {\bf x}_j ). 
\label{efficient_quadratic_reconstruction_final_wL}
\end{eqnarray} 
This is a highly efficient quadratic {\color{black} interpolation} formula because it only requires the solution at the cell and the gradients at the nodes of the cell. There is no need to explicitly compute and store second derivatives as typically required in a third-order finite-volume discretization on unstructured grids. Similarly, we have the following formula for the right state: 
\begin{eqnarray}
w_R  = w_k + \frac{1}{2}    \left(  {\bf g} + \overline{\bf g}_k  \right)  \cdot   ( {\bf x} - {\bf x}_k ),
\label{efficient_quadratic_reconstruction_final_wR}
\end{eqnarray}
where $ \overline{\bf g}_k$ is the arithmetic average of the gradients stored at nodes of the neighbor cell $k$. In the rest of the paper, these efficient quadratic {\color{black} interpolation} schemes will be referred to as {\color{black} NGQI (nodal-gradient quadratic interpolation)} schemes. 

Using the {\color{black} NGQI} scheme, we interpolate the left and right states at the face centroid $ {\bf x}_T$: 
\begin{eqnarray}
w_L &=&w_j +  \frac{1}{2}    \left(    \overline{\bf g}_T  +  \overline{\bf g}_j   \right)  \cdot   ( {\bf x}_T - {\bf x}_j ), 
\label{efficient_quadratic_reconstruction_final_wL_xT}  \\ [1ex]
w_R &=& w_k +  \frac{1}{2}    \left(   \overline{\bf g}_T  +  \overline{\bf g}_k   \right)  \cdot   ( {\bf x}_T - {\bf x}_k ), 
\label{efficient_quadratic_reconstruction_final_wR_xT}
\end{eqnarray} 
 where $ \overline{\bf g}_T  =  (  {\bf g}_{1} +  {\bf g}_{2} +   {\bf g}_{4} ) / 3$ for the shaded face in Figure \ref{tetra_quad_reconstruction_j_k}, and compute the numerical flux with them as in Eq. (\ref{numerical_flux_Phijk}). The flux Jacobian at a node $i$ as required in Eq. (\ref{flux_quadrature_CQ}) is computed with the arithmetic average of the interpolated solutions: 
\begin{eqnarray}
\frac{ w_{i,L} + w_{i,R} }{2},
\end{eqnarray} 
where in this case, the left and right states are computed at the node $i$,
\begin{eqnarray}
w_{i,L} &=&w_j +  \frac{1}{2}    \left(   {\bf g}_i +  \overline{\bf g}_j   \right)  \cdot   ( {\bf x}_i - {\bf x}_j ), \\ [1ex]
w_{i,R} &=& w_k +  \frac{1}{2}    \left(   {\bf g}_i +  \overline{\bf g}_k   \right)  \cdot   ( {\bf x}_i - {\bf x}_k ).
\label{efficient_quadratic_reconstruction_final_wLwR_node}
\end{eqnarray}
Here, the nodal gradient ${\bf g}_i$ is directly used and no interpolation is necessary.
 
Finally, we remark that the {\color{black} NGQI} schemes, Eqs. (\ref{efficient_quadratic_reconstruction_final_wL_xT}), and (\ref{efficient_quadratic_reconstruction_final_wR_xT}) at the face centroid will reduce to the FANG linear {\color{black} interpolation} \cite{NishikawaWhite_FANG:jcp2020,Nishikawa_FANG_AQ:Aviation2020,WhiteNishikawaBaurle_scitech2020,NishikawaWhite:AIAA2021-2720} if $\overline{\bf g}_j$ and $ \overline{\bf g}_k$ are replaced by the face averaged gradient $ \overline{\bf g}_T$: 
\begin{eqnarray}
w_L &=& w_j +    \overline{\bf g}_T   \cdot   ( {\bf x}_T - {\bf x}_j ),
\label{FANG_wL}
\\ [1ex]
w_R &=& w_k +     \overline{\bf g}_T    \cdot   ( {\bf x}_T - {\bf x}_k ),
\label{FANG_linear_reconstruction_final_wLwR_face_center}
\end{eqnarray}
or the cell-averaged nodal-gradient (CANG) linear {\color{black} interpolation} (as originally proposed by Zhang \cite{zhang_etal:2015,Zhang:Preprint2017,zhang:phd}) if $ \overline{\bf g}_T$ is replaced by $\overline{\bf g}_j$ in $j$ and $ \overline{\bf g}_k$ in $k$: 
\begin{eqnarray}
w_L &=& w_j +    \overline{\bf g}_j   \cdot   ( {\bf x}_T - {\bf x}_j ),
\label{CANG_wL}
\\ [1ex]
w_R &=& w_k +     \overline{\bf g}_k    \cdot   ( {\bf x}_T - {\bf x}_k ),
\label{CANG_linear_reconstruction_final_wLwR_face_center}
\end{eqnarray}
which show how simple it is to implement the proposed quadratic {\color{black} interpolation} scheme in an existing code if FANG or CANG are already available, and also how simple it is to revert to a second-order discretization. Furthermore, quite interestingly, we see that the proposed efficient quadratic {\color{black} interpolation} scheme is just the arithmetic average of the FANG and CANG linear {\color{black} interpolation} schemes: e.g., Eq. (\ref{efficient_quadratic_reconstruction_final_wL_xT}) is the average of Eqs (\ref{FANG_wL}) and (\ref{CANG_wL}). But of course, for the {\color{black} interpolation} scheme to be truly quadratically exact, the nodal gradients need to be computed by the quadratic LSQ method. 
\newline
\newline
{\bf Remark}: In the NASA VULCAN-CFD code \cite{WhiteNishikawaBaurle_scitech2020,VULCAN_FANGplus_scitech2022}, the FANG method is employed in combination with the linearity-preserving UMUSCL scheme \cite{nishikawa_LP_UMUSCL:JCP2020}: 
\begin{eqnarray}
w_L &=& \kappa \frac{  w_j  +  w'_k }{2} + ( 1 - \kappa) \left[ w_j +    \overline{\bf g}_T   \cdot   ( {\bf x}_T - {\bf x}_j ) \right], 
\label{FANG_wL_LP_UMUSCL}
\\ [1ex]
w_R &=&  \kappa \frac{  w'_j  +  w_k  }{2} + ( 1 - \kappa) \left[  w_k +     \overline{\bf g}_T    \cdot   ( {\bf x}_T - {\bf x}_k ) \right], 
\label{FANG_wR_LP_UMUSCL}
\end{eqnarray}
where
\begin{eqnarray}
w'_k  &=& w_k +   \overline{\bf g}_T  \cdot ( {\bf x}_T  + ( {\bf x}_T   - {\bf x}_j) -  {\bf x}_k ),
\\ [1ex] 
 w'_j &=& w_j +  \overline{\bf g}_T \cdot ( {\bf x}_T   + ( {\bf x}_T   - {\bf x}_k) -  {\bf x}_j ),
\end{eqnarray}
with $\kappa=1/3$. For comparison, we will employ this scheme as the baseline second-order scheme.


\subsection{Weak Boundary Conditions}

\begin{figure}[t] 
\begin{center}
\begin{minipage}[b]{0.8\textwidth}
\begin{center}
          \includegraphics[width=0.55\textwidth,trim=0 0 0 0,clip]{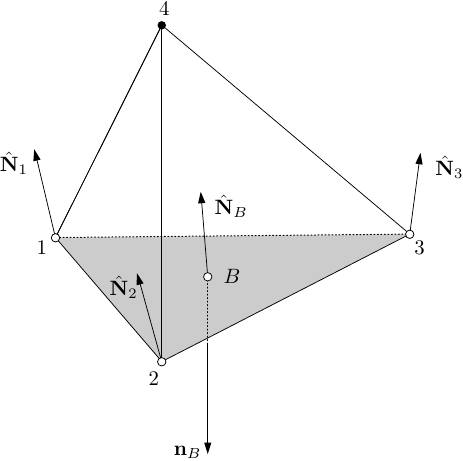}
        \caption{Tetrahedral cell adjacent a boundary, where the shaded triangular face is a boundary face, $ {\bf n}_B$ is an outward 
        normal vector of the boundary face, $ \hat{\bf N}_1$, $ \hat{\bf N}_2$, and $ \hat{\bf N}_3$ are unit vectors normal to the true boundary surface at boundary nodes $1$, $2$, and $3$, respectively, and $\hat{\bf N}_B =  \overline{\bf N} / |\overline{\bf N}| $, $ \overline{\bf N} = (  \hat{\bf N}_1 + \hat{\bf N}_2 + \hat{\bf N}_3 )/3  $.}
       \label{boundary_face_nodal_normals}
\end{center}
\end{minipage}
\end{center}
\end{figure}

For boundary conditions, we consider imposing them weakly through the right state in the numerical flux at a boundary face \cite{carlson:NASATM2011,mengaldo_etal_bcs:AIAA2014}:  
\begin{eqnarray}
{\bf w}_R  =   (   \rho_\infty , {\bf v}_\infty,   p_\infty), 
\end{eqnarray} 
for a free stream condition, where subscript $\infty$ indicates a free stream value, and 
\begin{eqnarray}
{\bf w}_R  =   (   \rho_L , {\bf v}_L,   p_\infty), 
\end{eqnarray} 
for an outflow condition (known as the back pressure condition) where $\rho_L$ and ${\bf v}_L$ are the values interpolated at the center of the boundary face from the cell adjacent to the boundary face using either Eq. (\ref{efficient_quadratic_reconstruction_final_wL_xT}) or Eq. (\ref{FANG_wL_LP_UMUSCL}). Note that it suffices to incorporate boundary conditions in the right state entering the numerical flux at the boundary face center. The flux correction term in Eq. (\ref{final_third_order_res}) is purely a mathematical term that eliminates a leading second-order error in the flux quadrature, and therefore, we simply set $w_{i,R} =w_{i,L}$ at a boundary node ($i = 1, 2, 3$ in Figure \ref{boundary_face_nodal_normals}) in evaluating the flux Jacobian $\left( {\partial {\bf f}}/{\partial {\bf w}} \right)_{i}$. 

It is known that a high-order curved grid is not necessary to achieve high-order accuracy but high-order accurate surface normal vectors are needed at a quadrature point  \cite{krivodonova_berger:JCP2006,nishikawa_boundary_quadrature:JCP2015}. For example, we implement a slip-wall condition relevant to inviscid problems considered in this paper by defining the right state as
\begin{eqnarray}
{\bf w}_R  =   \left(    \rho_L, {\bf v}_L -  2 (   {\bf v}_L  \cdot \hat{\bf N}_B   )  \hat{\bf N}_B  ,   p_L  \right), 
\label{slip_wall_wR}
\end{eqnarray}
where $\hat{\bf N}_B$ is a unit wall normal vector. It may not be immediately obvious, but this algorithm is equivalent to Algorithm II in Ref.~\cite{krivodonova_berger:JCP2006}. It is easy to show that the velocity in the direction $\hat{\bf N}_B$ vanishes on average at the boundary face centroid, which is the quadrature point in our discretization: 
\begin{eqnarray}
\frac{  {\bf v}_L +  {\bf v}_R } {2}  \cdot  \hat{\bf N}_B
&=&
\frac{  {\bf v}_L +  {\bf v}_L -  2 (   {\bf v}_L  \cdot \hat{\bf N}_B   )  \hat{\bf N}_B} {2}  \cdot  \hat{\bf N}_B
\nonumber \\ [2ex]
&=&
   {\bf v}_L  \cdot \hat{\bf N}_B  -  (   {\bf v}_L  \cdot \hat{\bf N}_B   ) \hat{\bf N}_B  \cdot \hat{\bf N}_B
\nonumber \\ [2ex]
&=&
 {\bf v}_L  \cdot \hat{\bf N}_B -  (   {\bf v}_L  \cdot\hat{\bf N}_B   ) 
\nonumber \\ [2ex]
&=&
  0.
\end{eqnarray}
Here, the unit normal vector $\hat{\bf N}_B$ needs to be exact for a quadratic geometry for the discretization to be third-order accurate \cite{krivodonova_berger:JCP2006,nishikawa_boundary_quadrature:JCP2015}, but the face centroid is not located exactly at a true surface, in general, and thus it is not exactly clear how to compute an accurate normal vector there. To avoid the ambiguity, we propose the following approach: first, compute accurate normal vectors at boundary nodes, which are typically defined at a true boundary surface, and then compute $\hat{\bf N}_B$ by a linear interpolation: 
\begin{eqnarray}
 \hat{\bf N} = \frac{ \hat{\bf N}_1 + \hat{\bf N}_2 + \hat{\bf N}_3   }{3},
 \label{averaged_boundary_normals}
\end{eqnarray}
where $\hat{\bf N}_1$, $\hat{\bf N}_2$, and $ \hat{\bf N}_3$ are sufficiently accurate unit surface normal vectors at boundary nodes $1$, $2$, and $3$, respectively, followed by the normalization: 
\begin{eqnarray}
 \hat{\bf N}_B = \frac{  \hat{\bf N} }{ | \hat{\bf N}|}.
 \label{averaged_boundary_normals_unit}
\end{eqnarray}
Linear interpolation is sufficient for quadratic exactness because normal vectors are derivatives of an algebraic equation defining a boundary surface (which should be quadratically exact, at least). In this work, we consider geometries for which exact surface normal vectors can be easily computed at each boundary node.

It is important to note that the interpolated surface normal vector $\hat{\bf N}_B$ is used only in the definition of the right state and the face normal vector ${\bf n}_B$ must be used otherwise, so that the residual, Eq. (\ref{final_third_order_res}), is properly closed and thus vanishes for a free stream flow, for example. It is emphasized again that a high-order curved cell is not required to achieve third-order accuracy and this is a well-known fact that has been demonstrated in other high-order methods \cite{krivodonova_berger:JCP2006,nishikawa_boundary_quadrature:JCP2015}. See especially Ref.\cite{krivodonova_berger:JCP2006}, which considers a similar high-order finite-volume-type method and demonstrates up to fifth-order accuracy with linear grids for problems with a curved boundary. 


 {\color{black}
\subsection{Comments on Comparison with Other Third-Order Finite-Volume Methods}

It is beyond the scope of the paper to make a comprehensive comparison with standard third-order cell-centered finite-volume methods, but several remarks are in order. First, the main advantage of the proposed third-order method is that there is no need to compute and store second derivatives, which is required in standard third-order cell-centered finite-volume methods. It means one can avoid the computation and storage of 30 second derivatives (six for each variable) in the case of the Euler equations in three dimensions. The saving here is by a factor of three: 15 (gradient components for the five variables) in our method and 45 (15 gradient components and 30 second derivatives) in standard methods. 

Second, the use of nodal gradients is already more efficient than standard methods, where gradients and Hessians are computed and stored at cells, for tetrahedral grids approximately by a factor of six (two in two dimensions) \cite{NishikawaWhite_FANG:jcp2020,WhiteNishikawaBaurle_scitech2020}. Therefore, in total, the proposed third-order method can be more efficient by a factor of 18 than standard methods, in the gradient and Hessian computations and storage. 

Third, the quadratic interpolation formula is more efficient as it does not require the computation of a curvature term involving a Hessian matrix, e.g., $ ( {\bf x}_T - {\bf x}_j )^t {\nabla}^2 w_j   ( {\bf x}_T - {\bf x}_j )$. Hence, we avoid a matrix-vector product followed by a dot product of two vectors, i.e., 60 multiplications and 40 additions, in total, in the quadratic interpolation for five solution variables. 


Finally, standard methods typically employ Gaussian quadrature formulas for the flux integration, requiring at least three numerical flux evaluations per face \cite{RokickiWieteska_ECCOMAS2006} or six \cite{ZhuShu:JCP2020}, whereas we only need to compute a single numerical flux per face. The factor of three can be significant for an expensive numerical flux function as the flux correction term can be much cheaper to compute. The same can be true for high-order source integrations. Certainly, it would be possible to derive similar integration formulas in the form of a point value with a correction also in standard methods, but the correction will involve second derivatives \cite{PintBrennerCinnellaMaugarsRobinet:JCP2017,SetzweinEssGerlinger:JCP2021}. 

In this paper, we do not present any comparative study because the whole point of this work is to avoid implementing a standard third-order method that requires second derivatives and multiple quadrature points. Such a study should be performed carefully with an established code of a standard third-order method (clearly defining what `standard' means) and is left as a future work.
}

{\color{black}

\section{Nonlinear Solver and Time Integration Scheme}
\label{Solver_TimeScheme}

\subsection{Nonlinear Solver}
\label{NonlinearSolver}

In this work, for steady problems, we solve a global system of nonlinear residual equations,
\begin{eqnarray}
0 = {\bf Res} ( {\bf U} ),
\label{res_eq_ns}
\end{eqnarray}
where ${\bf Res} $ denotes a global vector of the spatial part of the discretization (\ref{final_third_order_res}) and $ {\bf U} $ denotes a global vector of the conservative variables, by the implicit defect-correction solver:
\begin{eqnarray}
{\bf U}^{n+1}  = {\bf U}^{n}  +  \Delta {\bf U},
\end{eqnarray}
\begin{eqnarray}
 \frac{\partial  \overline{   {\bf {Res} } }  }{ \partial {\bf U}} \Delta {\bf U} = -  {\bf Res} ( {\bf U}^{n} ), 
\label{linearized_system}
\end{eqnarray}
where $n$ is the iteration counter and the Jacobian $\partial \overline{{ \bf {Res} } } / \partial {\bf U}$ is the exact differentiation of 
the low-order compact residual $\overline{\bf Res}$ with zero nodal gradients (i.e., first-order accurate residual). The linear system is relaxed by the block multi-color Gauss-Seidel relaxation scheme, which is performed in the order of colors with the following scheme for a cell $j$,
\begin{eqnarray}
 \Delta {\bf u}^{m+1}_j  =  \Delta {\bf u}^m_j + \omega \,  {\bf r}_j,  
 \quad 
 {\bf r}_j = \left( \frac{\partial \overline{{ \bf {Res} }_j }}{ \partial {\bf u}_j} \, \right)^{-1} 
 \left[  - \sum_{k \in \{ k_j \}   }   \frac{\partial \overline{{ \bf {Res} }_j } }{ \partial {\bf u}_k} \Delta {\bf u}^m_k  - {\bf Res}_j ( {\bf U}^n )  \right] -  \Delta {\bf u}^m_j ,
\label{idc:GS_flow}  
\end{eqnarray}
where the set $\{ k_j \}$ contains all the face-neighbor cells, $ \omega=1$, ${ \bf {Res} }_j$ is the residual vector in the cell $j$, $m$ is the relaxation counter, and $ \frac{\partial \overline{{ \bf {Res} }_j } }{ \partial {\bf u}_j} $ and $ \frac{\partial \overline{{ \bf {Res} }_j } }{ \partial {\bf u}_k} $ are the 5$\times$5 diagonal and off-diagonal matrices, respectively. In this paper, we will perform only five relaxations per nonlinear iteration ($m=5$) unless otherwise stated.

\subsection{Time Integration Scheme}
\label{TimeIntegration}

For unsteady problems, we employ the explicit third-order SSP-RK time integration scheme \cite{SSP:SIAMReview2001}:
\begin{eqnarray}
{\bf U}^{(1)}  &=& {\bf U}^{n}  -  \Delta  t \,   \left\{   {\bf M}^{-1}  {\bf Res}' (  {\bf U}^n,  t^n   ) -  {\bf S}(  t^n )   \right\} ,      \\ [2ex]
{\bf U}^{(2)}  &=&   \frac{ 3 {\bf U}^{n}  +  {\bf U}^{(1)}  }{4}  -  \frac{ \Delta t }{4}  \left\{    {\bf M}^{-1}  {\bf Res}' ( {\bf U}^{(1)},  t^n + \Delta  t )  -  {\bf S}(  t^n + \Delta  t) \right\},    \\ [2ex]
{\bf U}^{n+1}  &=& \frac{    {\bf U}^{n}  + 2  {\bf U}^{(2)}  }{3}  -   \frac{2  \Delta t }{3}   \left\{  {\bf M}^{-1}  {\bf Res}' ( {\bf U}^{(2)},  t^n + \Delta  t /2 )  -  {\bf S}(  t^n + \Delta  t/2)  \right\},     
\end{eqnarray} 
where $t^n$ denotes the physical time at the $n$-th time level, $n=0,1,2, \cdots$, $ \Delta  t$ is a time step, ${\bf S}$ is a global vector of the source terms at cell centers, ${\bf S} = ( s_1, s_2, \cdots )$, and ${\bf Res}'$ denotes the flux-balance part of the spatial residual whose $j$-th component is given by
\begin{eqnarray}
{\bf Res}'_j ( {\bf U}^n,  t^n  )
=  
\frac{1}{V_j}
 \sum_{ k \in \{  k_j \}  }  
   \left[  
 \Phi_{jk}
  +
  \frac{1}{24}     \sum_{i=1}^3   \left(  \frac{\partial {\bf f}}{\partial {\bf w}} \right)_{\!\! i}  \!\!    \nabla {\bf w}_i   ( {\bf x}_i - {\bf x}_{T} ) \right]  | {\bf n}_T |,
\end{eqnarray} 
where the boundary condition and the source terms are evaluated at $t=t^n$. 
The symbol $ {\bf M}^{-1} $ denotes the inversion of the mass matrix arising from the source quadrature formula applied to the time derivative terms: e.g., the $j$-th component of the product of the mass matrix and a global vector of arbitrary function values at cell centers, ${\bf z} = ( z_1, z_2, \cdots )$ is given by 
\begin{eqnarray}
\left(  {\bf M} {\bf z} \right)_j
=
 {\bf z}_j  V_j 
  + 
   \frac{1}{40}
   \sum_{i=1}^4
   \nabla {\bf z}_i     ( {\bf x}_i - {\bf x}_{j} ),
\end{eqnarray} 
where the nodal gradient $ \nabla {\bf z}_i $ is computed by the quadratic LSQ method as described in Section \ref{LSQ_stencils}. The inversion can be efficiently performed for an arbitrary vector ${\bf b}$ as ${\bf M}^{-1} {\bf b}$, without storing the matrix ${\bf M}$, by the fixed-point iteration: for $\ell = 0, 1, 2, \cdots$.
\begin{eqnarray}
{\bf z}^{\ell+1} = {\bf z}^{\ell}  -  {\bf V}^{-1} {\bf r},  
\end{eqnarray} 
with the initial value
\begin{eqnarray}
{\bf z}^{0} =   {\bf V}^{-1}  {\bf b},
\end{eqnarray} 
where ${\bf V}$ is a diagonal matrix of cell volumes $V_j$ and ${\bf r} = {\bf M} {\bf z}^{\ell} -  {\bf b}  $ denotes the fixed-point iteration residual: the $j$-th component is given by
\begin{eqnarray}
{\bf r}_j = 
 {\bf z}_j^{\ell}  V_j
  + 
   \frac{1}{40}
   \sum_{i=1}^4
  \nabla {\bf z}_i^{\ell}  
   ( {\bf x}_i - {\bf x}_{j} )  V_j
   - 
 {\bf b}_j.
\end{eqnarray}
At convergence, we will have ${\bf z} = {\bf M}^{-1} {\bf b}$. 
In this study, this iteration method is used to compute ${\bf M}^{-1} {\bf Res}'$ in the explicit third-order SSP-RK time integration scheme. Note that the iteration converges at the first iteration for the point evaluation of the source term (i.e., ${\bf M} = {\bf V}$ and thus the source correction is zero) and simply gives the division by the cell volume. For our high-order source quadrature formula, this iteration scheme converges very rapidly. For all the problems considered in this paper, it was sufficient to perform two or three fixed-point iterations (which reduces the fixed-point iteration residual by two or three orders) to observe third-order accuracy. Therefore, it only requires two or three evaluations of the high-order source quadrature.

}
 
\section{Results}
\label{Results}

In this paper, we focus on simple but illustrative test cases with relatively coarse grids due to the serial nature of the testbed code. Steady problems are solved by the nonlinear solver described in Section \ref{NonlinearSolver} and unsteady problems are solved by the explicit time-integration scheme described in Section \ref{TimeIntegration}. For all the test cases, we visualize the results based on the solution values recovered at nodes by averaging the values interpolated from the cells sharing a node. For example, for a scalar variable $w$ at a node $i$, the $\kappa=0$ {\color{black} interpolation} scheme is used for the baseline scheme, 
 \begin{eqnarray}
{w}_i  = 
\frac{1}{ n_i  }
\sum_{k \in    \{   k_i^{(0)}  \} } 
\left[ {w}_k +  \overline{\bf g}_k    \cdot   ( {\bf x}_i - {\bf x}_k ) \right],
\end{eqnarray}
where, as defined earlier, $ \{   k_i^{(0)}  \} $ is the set of cells sharing the node $i$, $n_i$ denotes the number of cells in the set, ${ w}_k$ is the solution stored at the cell $k$, and $\kappa=0$ is chosen here to avoid the ambiguity of selecting a cell on the other side, e.g., $k$ in $w_L$, Eq. (\ref{FANG_wL_LP_UMUSCL}), and the {\color{black} NGQI} scheme for others,
 \begin{eqnarray}
{w}_i  = 
\frac{1}{ n_i }
\sum_{k \in    \{   k_i^{(0)}  \} } 
\left[ {w}_k + \frac{1}{2}    \left(  {\bf g}_i + \overline{\bf g}_k  \right)  \cdot   ( {\bf x}_i - {\bf x}_k ) \right].
\end{eqnarray}
Therefore, the design order of accuracy is maintained in the solutions recovered at nodes for all cases.

\newpage

\subsection{Accuracy Verification} 

\subsubsection{Dirichlet condition in rectangular domain (three-dimensional solution)} 

\begin{figure}[htbp!] 
  \begin{center}
    \centering
          \begin{subfigure}[t]{0.48\textwidth}        
        \includegraphics[width=0.85\textwidth]{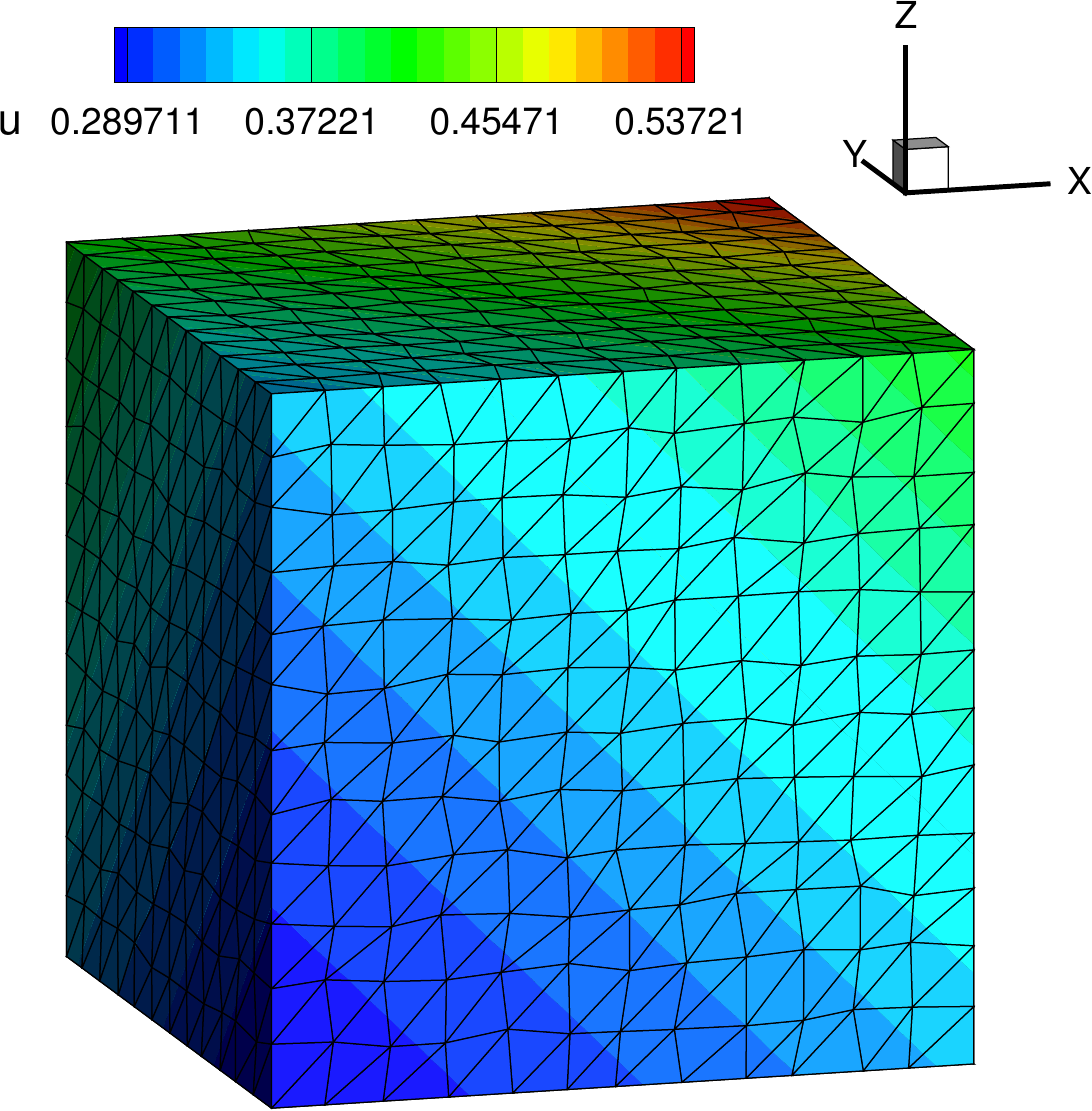}
          \caption{Coarsest grid.}
          \label{fig:accv_rectangular_exp_grid01}
      \end{subfigure}
      \hfill
          \begin{subfigure}[t]{0.48\textwidth}        
        \includegraphics[width=0.99\textwidth]{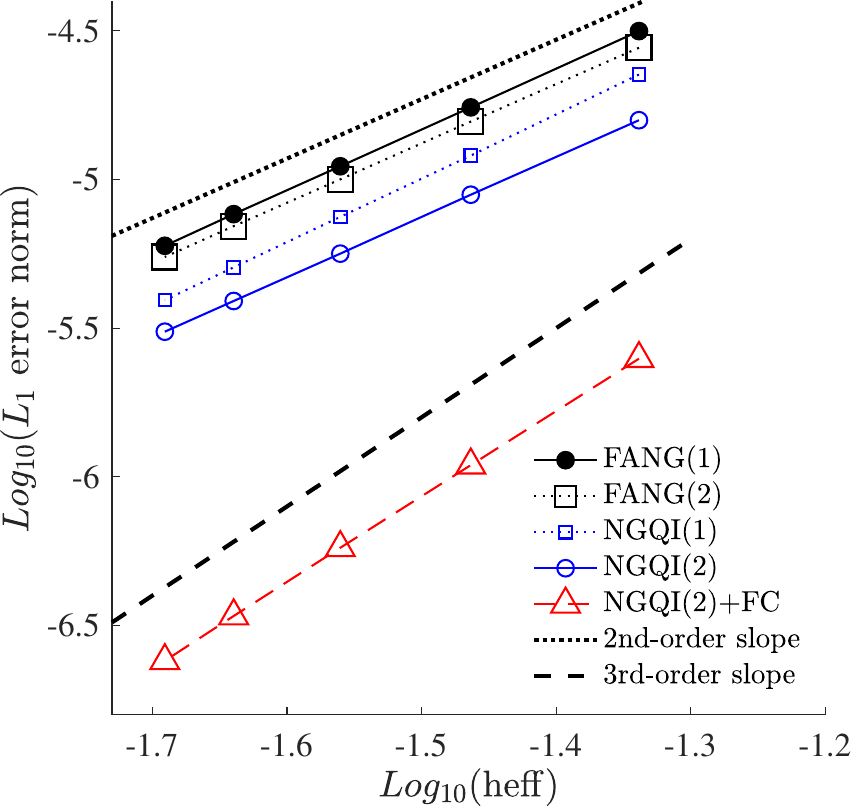}
          \caption{Error convergence.}
          \label{fig:accv_rectangular_exp_err_conv}
      \end{subfigure}
      \hfill
        \caption{Accuracy verification result for a cube domain: (a) the coarsest grid with contours of the $x$-velocity, (b) error convergence results for the $x$-velocity.
        }
          \label{fig:accv_rectangular_exp}
\end{center} 
\end{figure}

{\color{black} We first perform an accuracy verification study by the method of manufactured solutions \cite{Roy_etal_IJNMF2007} using the following three-dimensional solution set: 
\begin{eqnarray}
\rho = 1 + \exp \left( s_{xyz} \right), \quad
    p = 1 + \exp \left(s_{xyz}  \right), \quad
s_{xyz} =   0.5 ( x+y+z) , 
\\
u  = 0.2 + \exp \left(  s_{xyz}  \right), \quad
v  = 0.2 + \exp \left(  s_{xyz}  \right), \quad
w  = 0.2 + \exp \left( s_{xyz}  \right), 
\end{eqnarray}
which is made an exact solution to the compressible Euler equations by introducing source terms that are numerically evaluated at the cell centroid and added, together with the source correction, to the residual at each cell. The derivatives of the source terms were computed numerically by repeated uses of the chain rule from the derivatives of the solutions above. Here, we consider a unit cube domain as shown in Figure \ref{fig:accv_rectangular_exp_grid01} and solve the Euler equations with a weak Dirichlet boundary condition, where the right state at a boundary face is given by the exact solution evaluated at the face centroid, over a series of relatively coarse irregular tetrahedral grids with $n^3$ nodes, $n=13, 17, 21, 25, 28$, having 10,368, 24,576, 48,000, 82,944, and 118,098 tetrahedra, respectively. Five different discretizations are considered. FANG(1) and FANG(2) refer to the baseline FANG schemes using the LP-UMUSCL schemes, Eqs. (\ref{FANG_wL_LP_UMUSCL}), and (\ref{FANG_wR_LP_UMUSCL}) with linear and quadratic LSQ gradients, respectively. The {\color{black} NGQI}(1) and {\color{black} NGQI}(2) refer to the {\color{black} NGQI} schemes with linear and quadratic LSQ gradients, respectively. All these schemes are second-order accurate with the one-point flux quadrature formula as in Eq. (\ref{2nd-order_CCFV}). Finally, {\color{black} NGQI}(2)+FC refers to the proposed economical third-order scheme, which adds the flux correction as in Eq. (\ref{flux_quadrature_CQ}) to {\color{black} NGQI}(2). For the quadratic LSQ method, the LSQ stencil is augmented with node-neighbors of the cells around each node, as mentioned earlier. For this problem, no augmentation was performed for the linear LSQ gradient method. The residual equations are solved by the implicit defect-correction solver to reduce the $L_1$ residual norms by ten orders of magnitude on each grid. As can be seen in Figure \ref{fig:accv_rectangular_exp_err_conv}, third-order accuracy is achieved only by NGQI(2)+FC, as expected. FANG(2) gives slightly lower errors, but remains second-order accurate. Similarly, NGQI(1) and NGQI(2) also produce lower errors but remain second-order accurate. 
}

\subsubsection{Slip condition over a curved boundary} 

\begin{figure}[htbp!] 
  \begin{center}
    \centering
          \begin{subfigure}[t]{0.48\textwidth}        
        \includegraphics[width=0.99\textwidth]{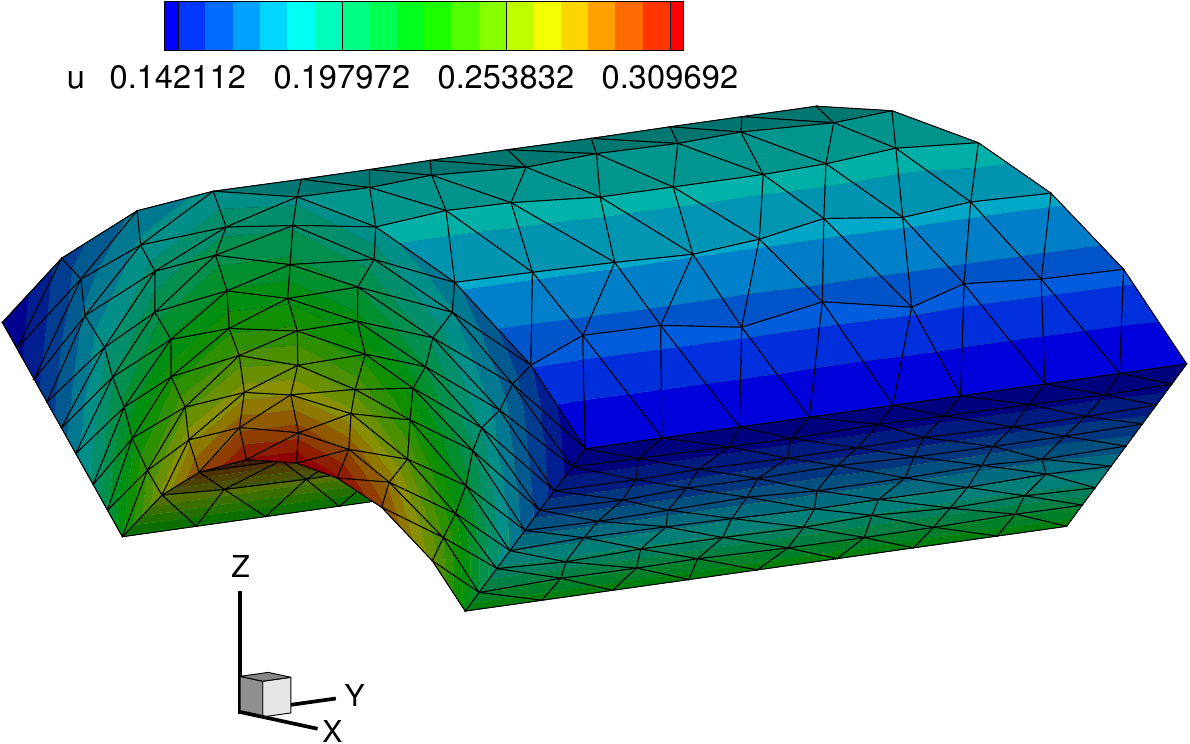}
          \caption{Coarsest grid.}
          \label{fig:accv_cyl_grid01}
      \end{subfigure}
      \hfill
          \begin{subfigure}[t]{0.48\textwidth}        
        \includegraphics[width=0.99\textwidth]{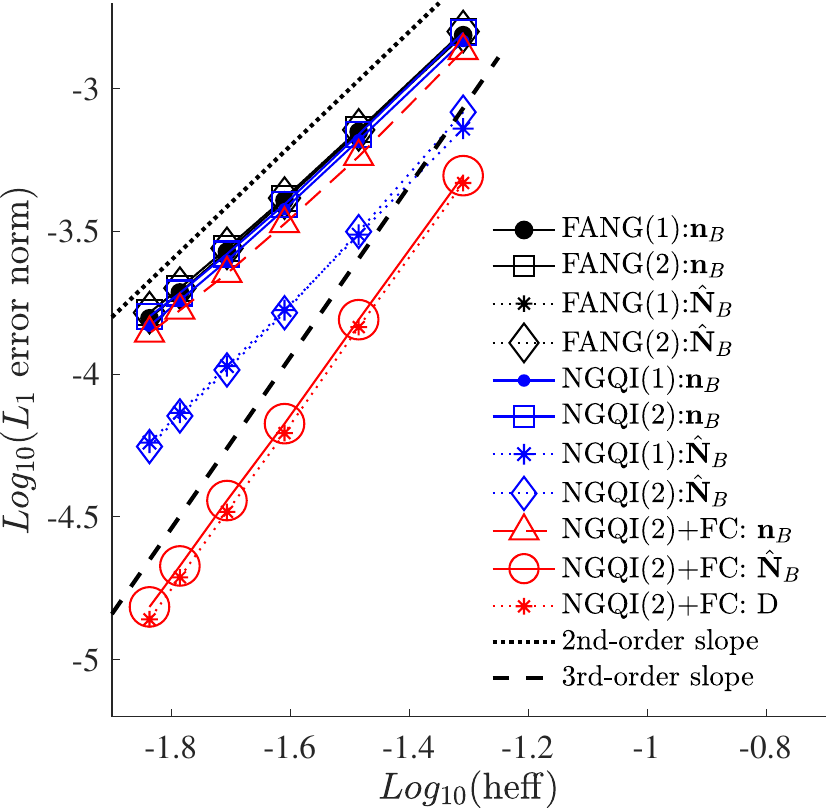}
          \caption{Error convergence.}
          \label{fig:accv_cyl_err_conv}
      \end{subfigure}
      \hfill 
        \caption{Accuracy verification result for a curved domain: (a) the coarsest grid with the $x$-velocity contours, (b) error convergence results for the density. The data labeled as {\color{black} NGQI}(2)+FC: D is the third-order scheme with the weak Dirichlet condition (the exact solution is given as the right state). 
        }
          \label{fig:accv_cyl}
\end{center} 
\end{figure}

{\color{black} To demonstrate third-order accuracy over a curved domain with a slip-wall condition, we consider a potential vortex solution:
\begin{eqnarray}
 u = \frac{ K z}{ x^2  +z^2 }, \quad
 v = 0, \quad
 w =   \frac{ - K x}{ x^2  +z^2 }, \quad
 p = \frac{ u^2 + w^2 }{2} + \frac{1}{\gamma}, \quad
 \rho = 1 + 0.1 \frac{ 1 - \exp( -1.5 ( r - 0.5 )  )  }{ 1  - \exp(- 1.5) }, 
\end{eqnarray}
where $r = \sqrt{ x^2 + y^2 + z^2} $ and $K=1 / ( 2 \pi) $ for a cylindrical domain, where the potential vortex solution is tangent to the bottom boundary. The coarsest grid is shown in Figure \ref{fig:accv_cyl_grid01}, with contours of the $x$-velocity. The same five schemes are tested as in the previous section: FANG(1), FANG(2), NGQI(1), NGQI(2), and NGQI(2)+FC}, over a series of irregular tetrahedral grids with $n^3$ nodes, $n=9, 13, 17, 21, 25, 28$, having 3,072, 10,368, 24,576, 48,000, 82,944, and 118,098 tetrahedra, respectively. To investigate the impact of the surface normal vector in the slip-wall boundary condition, we consider two choices for $\hat{\bf N}_B$ in Eq. (\ref{slip_wall_wR}): the boundary face normal ${\bf n}_B$ and the averaged exact surface normals over the nodes of the boundary face $\hat{\bf N}_B$ as defined by Eqs. (\ref{averaged_boundary_normals}) and (\ref{averaged_boundary_normals_unit}). At a boundary node $(x_i,y_i,z_i)$, the exact surface normal vector is given by $\hat{\bf N}_i =  {\bf N}_i / |{\bf N}_i|$, ${\bf N}_i = (x_i,0, z_i)$.

Error convergence results are shown in Figure \ref{fig:accv_cyl_err_conv} for the $x$-velocity component. Again, the errors are computed at the centroids against the exact solutions at the centroids (not cell averages). First of all, we observe no significant differences in the baseline FANG schemes with linear and quadratic LSQ gradients as well as with the two choices of the surface normal vectors. On the other hand, the {\color{black} NGQI} schemes show significantly lower errors with the more accurate surface normal $\hat{\bf N}_B$ than the boundary face normal ${\bf n}_B$. Finally, the third-order scheme, {\color{black} NGQI}(2)+FC, achieves third-order accuracy with the averaged nodal normals $\hat{\bf N}_B$; the errors are nearly at the same level as those obtained with the Dirichlet condition (using the exact solution in the right state), which is included for comparison and labeled as {\color{black} NGQI}(2)+FC: D in the legend. However, if used with the boundary face normal ${\bf n}_B$, it reduces to second-order accurate and the errors are close to those of the baseline second-order schemes. These results indicate that third-order accuracy can be achieved with linear grids as long as the surface normals are accurate (confirming the results of Krivodonova and Berger \cite{krivodonova_berger:JCP2006}) and that benefits of the quadratic solution {\color{black} interpolation} and the third-order scheme are lost if used with a lower-order surface normal.

\subsubsection{Unsteady problem with a manufactured solution (three-dimensional solution)} 

\begin{figure}[htbp!] 
  \begin{center}
    \centering
          \begin{subfigure}[t]{0.45\textwidth}        
        \includegraphics[width=0.85\textwidth]{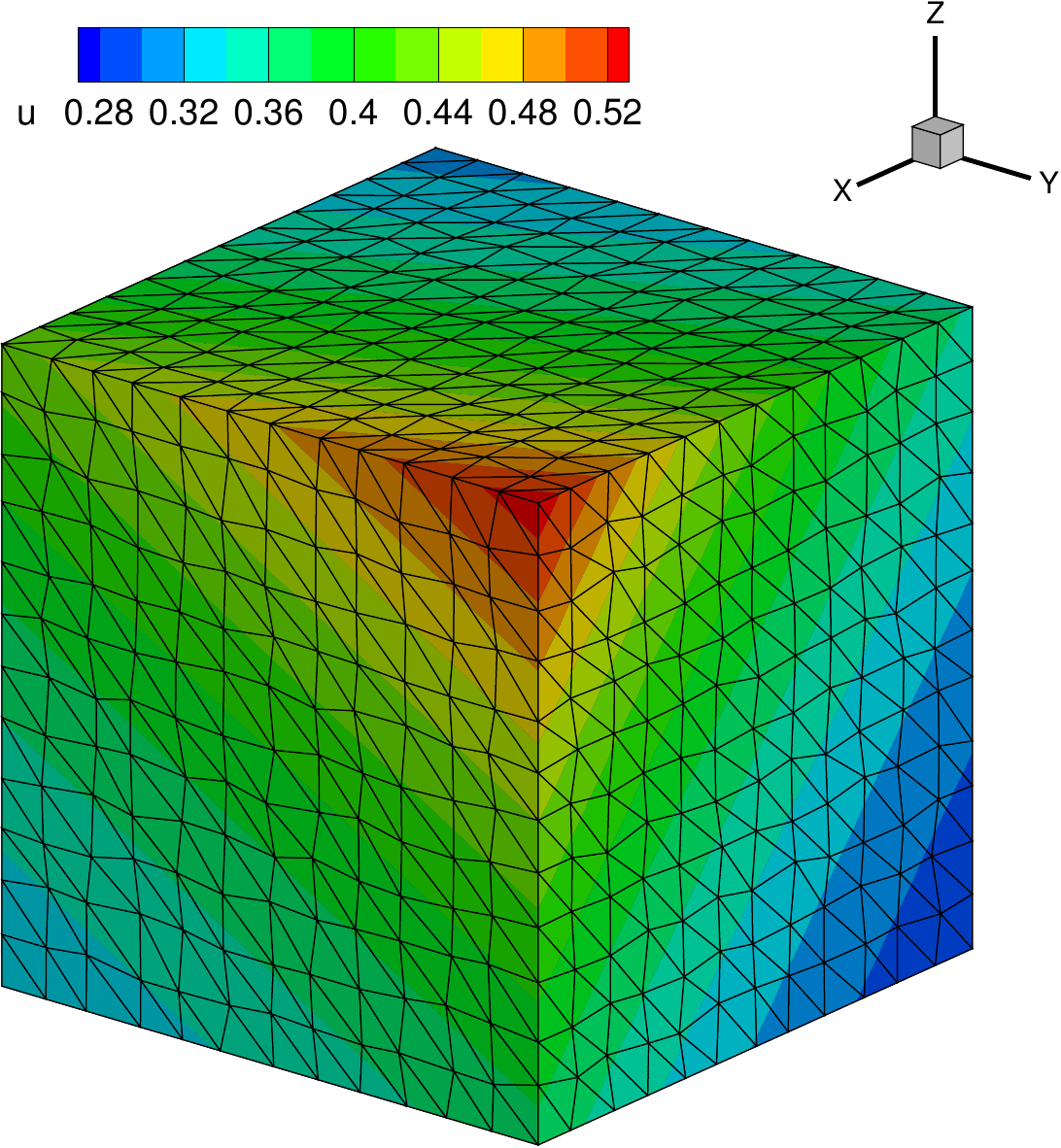}
          \caption{Coarsest grid.}
          \label{fig:accv_unsteady_mms_grid01}
      \end{subfigure}
      \hfill
          \begin{subfigure}[t]{0.45\textwidth}        
        \includegraphics[width=0.99\textwidth]{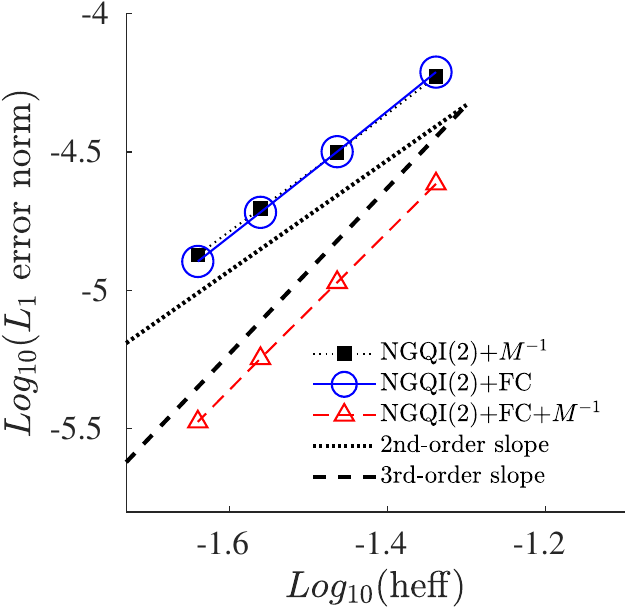}
          \caption{Error convergence.}
          \label{fig:accv_unsteady_mms_err_conv} 
      \end{subfigure}
      \hfill
        \caption{Accuracy verification result for a cube domain: (a) the coarsest grid with contours of the $x$-velocity, (b) error convergence results for the $x$-velocity.
        }
          \label{fig:accv_unsteady_mms}
\end{center} 
\end{figure}

{\color{black} To demonstrate third-order accuracy for unsteady problems, we employ the method of manufactured solutions \cite{Roy_etal_IJNMF2007} using the following three-dimensional solution set: 
\begin{eqnarray}
\rho = 1 + 0.1  \exp \left( s_{xyzt} \right), \quad
    p = 1.3  + 1.2 \exp \left(s_{xyzt}  \right), \quad  s_{xyzt} =   0.75 x + 0.35 y +   0.65  z  + 0.65 t, 
\\
u  = 0.2 + 0.1 \exp \left(  s_{xyzt}  \right), \quad
v  = 0.1 + 0.2 \exp \left(  s_{xyzt}  \right), \quad
w  = 0.3 + 0.3 \exp \left( s_{xyzt}  \right).
\end{eqnarray}
Starting from the solution at $t=0$ in a cubic domain, we perform the time integration up to the final time $t = 0.5$ by the third-order SSP-RK scheme as described in Section \ref{TimeIntegration} using $\Delta t = \, $CFL $ \min_{cells}  V_j  \{  \sum_{ k \in \{ k_j \} }  ( |u_n| + a )  | {\bf n}_T |/2  \}^{-1} $, where $u_n$ and $a$ are the face-normal velocity and the speed of sound evaluated at the Roe average state at the $k$-th face, with CFL $ = 0.95$. At boundaries, the weak Dirichlet condition is employed as in the previous tests. To measure the order of accuracy, we perform the computation over four levels of grids with 10,368, 24,576, 48,000, and 82,944 tetrahedra. The last time step is adjusted to finish the computation exactly at $t = 0.5$, which required 185, 248, 327, and 413 time steps in total for the four levels of grids, respectively. Figure \ref{fig:accv_unsteady_mms_grid01} shows the contours of the exact solution $u$ in the coarsest grid. 

To emphasize the importance of the flux correction and the mass-matrix inversion, we compare three schemes: NGQI(2)+$M^{-1}$, NGQI(2)+FC, and NGQI(2)+FC+$M^{-1}$, where NGQI(2)+$M^{-1}$ is the third-order scheme without the flux correction, NGQI(2)+FC is the third-order scheme without the mass-matrix inversion, and NGQI(2)+FC+$M^{-1}$ is the third-order scheme. As can be seen in the results shown in Figure \ref{fig:accv_unsteady_mms_err_conv}, third-order accuracy cannot be achieved if either the flux correction or the mass-matrix inversion is missing; third-order accuracy is achieved only when both are incorporated. 

}

\begin{figure}[htbp!] 
  \begin{center}
    \centering
          \begin{subfigure}[t]{0.45\textwidth}        
        \includegraphics[width=0.85\textwidth]{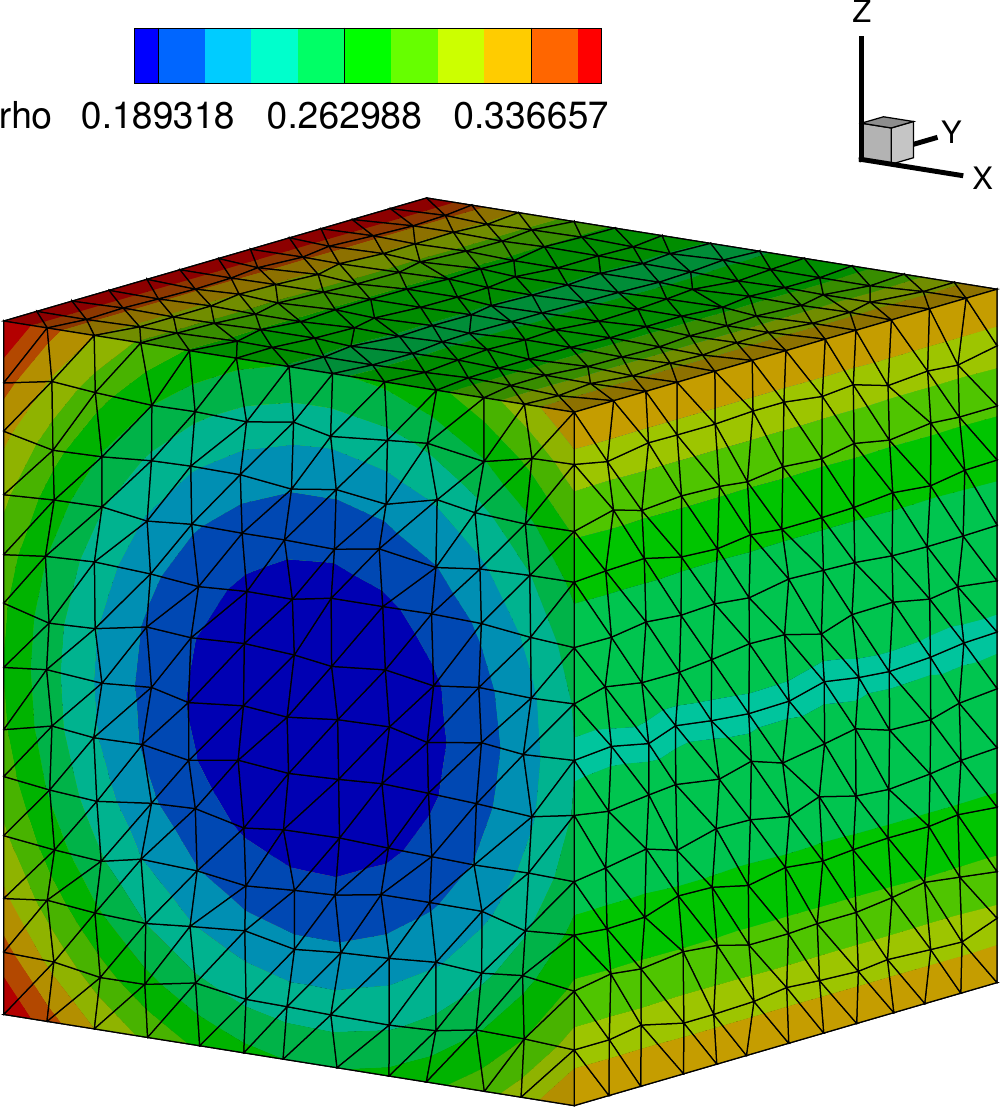}
          \caption{Coarsest grid.}
          \label{fig:accv_unsteady_vortex_grid01}
      \end{subfigure}
      \hfill
          \begin{subfigure}[t]{0.45\textwidth}        
        \includegraphics[width=0.99\textwidth]{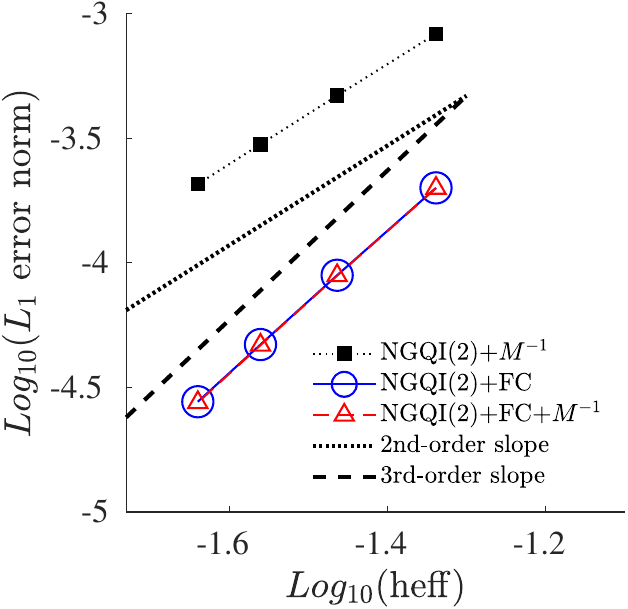}
          \caption{Error convergence.}
          \label{fig:accv_unsteady_vortex_err_conv}
      \end{subfigure}
      \hfill
        \caption{Accuracy verification result for a cube domain: (a) the coarsest grid with contours of the density, (b) error convergence results for the $x$-velocity.
        }
          \label{fig:accv_unsteady_vortex_exp}
\end{center} 
\end{figure}

\subsubsection{Unsteady problem with an exact inviscid vortex transport solution (two-dimensional solution)} 
\label{inviscid_vortex_accuracy_verification}

{\color{black} Finally, to show that the effect of the high-order source or equivalently the mass-matrix inversion on a two-dimensional problem, we consider an exact inviscid vortex transport solution to the Euler equations: 
\begin{eqnarray}
u =  u_\infty  -  \frac{ K \overline{y} }{2 \pi }   \exp \left(    \frac{ \alpha( 1-\overline{r}^2 ) }{2} \right) ,  \quad 
v = 0, \quad 
w  = w_\infty   + \frac{ K \overline{x} }{2 \pi }   \exp \left(      \frac{  \alpha( 1-\overline{r}^2 ) }{2} \right), 
\\ [1ex]
T =  1  -  \frac{ K^2 (\gamma-1)  }{8 \alpha \pi^2 }   \exp \left[   \alpha( 1-\overline{r}^2 )\right] , \quad
\rho = T^{  \frac{1}{\gamma-1}  } , \quad
p = \frac{ \rho ^{  \gamma   }  }{\gamma}    ,
\\ [1ex]
 \overline{x} = x  - u_\infty t, \quad
 \overline{z} = z  - w_\infty t, \quad
 \overline{r} ^2 =\overline{x}^2 + \overline{z}^2 ,
 \label{inviscid_vortex_sol}
\end{eqnarray}
where $(u_\infty,v_\infty,w_\infty)=(0.1,0.0,0.0)$, $\alpha=0.8$, and $K=6$. Note that a large parameter $K=6$ is chosen to avoid a false accuracy verification problem discussed in Ref.~\cite{Nishikawa_FakeAccuracy:2020}. This is an exact solution to the Euler equations in two dimensions, and therefore, there is no source term in this problem. The mass-matrix inversion is still necessary for the time derivative, but as we will show, it has no impact on the solution accuracy in this two-dimensional solution. To perform the accuracy verification efficiently, we restrict the domain to $(x,y,z) \in [ -0.5, 0.5] \times  [ -0.5, 0.5]   \times  [ -0.5, 0.5]   $, and perform the time integration with CFL $  = 0.95$. At boundaries, the weak Dirichlet condition is employed as in the previous tests. The last time step is adjusted to finish the computation exactly at the final time $t  = 0.5$, which required, in total, 80, 116, 142, 196 time steps over four levels of grids with 10,368, 24,576, 48,000, 82,944 tetrahedra, respectively. The exact density contours on the coarsest grid are shown in Figure \ref{fig:accv_unsteady_vortex_grid01}. As in the previous problem, we compare three schemes: NGQI(2)+$M^{-1}$, NGQI(2)+FC, and NGQI(2)+FC+$M^{-1}$. Figure \ref{fig:accv_unsteady_vortex_err_conv} shows the error convergence results. Clearly, in this case, the scheme without the mass-matrix inversion, NGQI(2)+FC, achieves third-order accuracy and produces nearly the same results as the third-order scheme, NGQI(2)+FC+$M^{-1}$. These results show that the flux correction is still required but the mass-matrix inversion has almost no impact on the solution accuracy for a two-dimensional problem. In fact, we have found that the source correction for the time derivative is negligibly small in this problem. This test case indicates that it is important to perform accuracy verification using a three-dimensional solution as in the previous section.
}

\subsection{Subsonic Flow over a Smooth Bump} 

\begin{figure}[htbp!] 
  \begin{center}
    \centering
          \begin{subfigure}[t]{0.48\textwidth}        
        \includegraphics[width=0.99\textwidth]{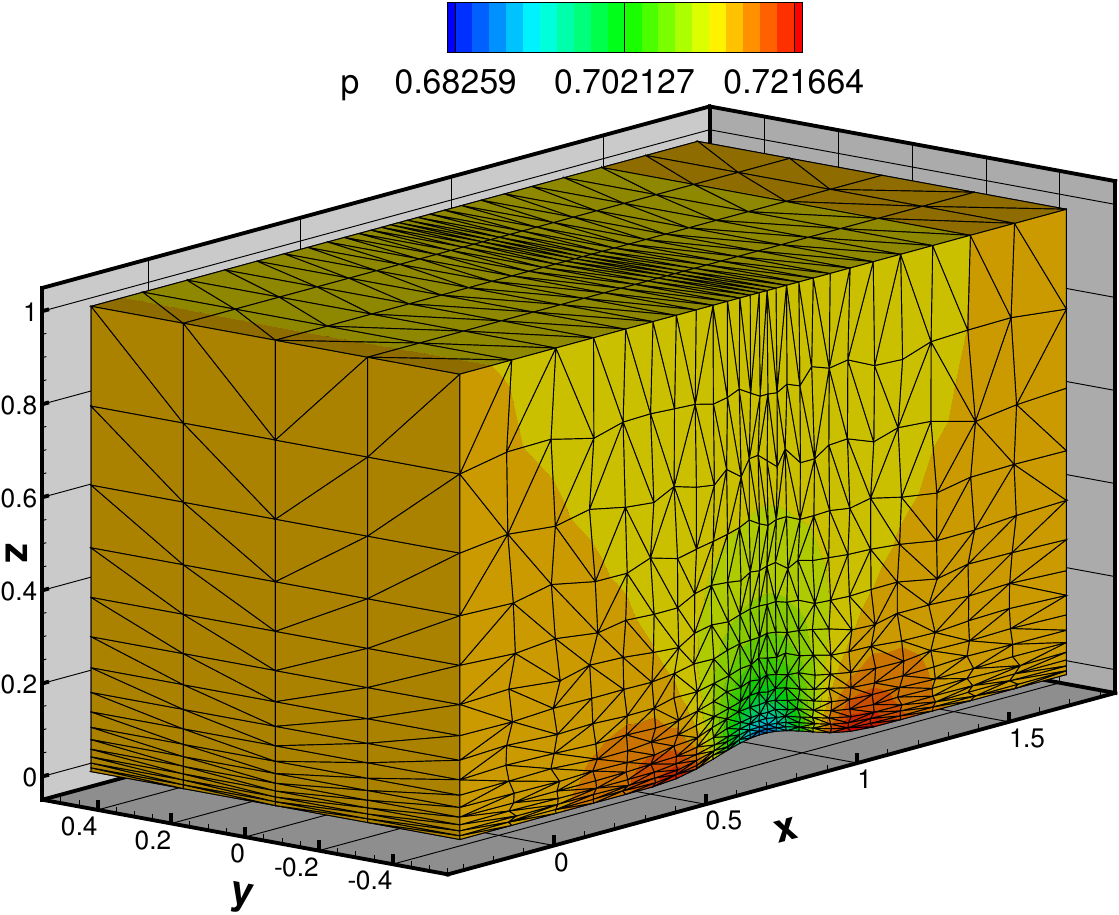}
          \caption{Coarsest grid with pressure contours.}
          \label{fig:bump_grid}
      \end{subfigure}
      \hfill
          \begin{subfigure}[t]{0.48\textwidth}        
        \includegraphics[width=0.99\textwidth]{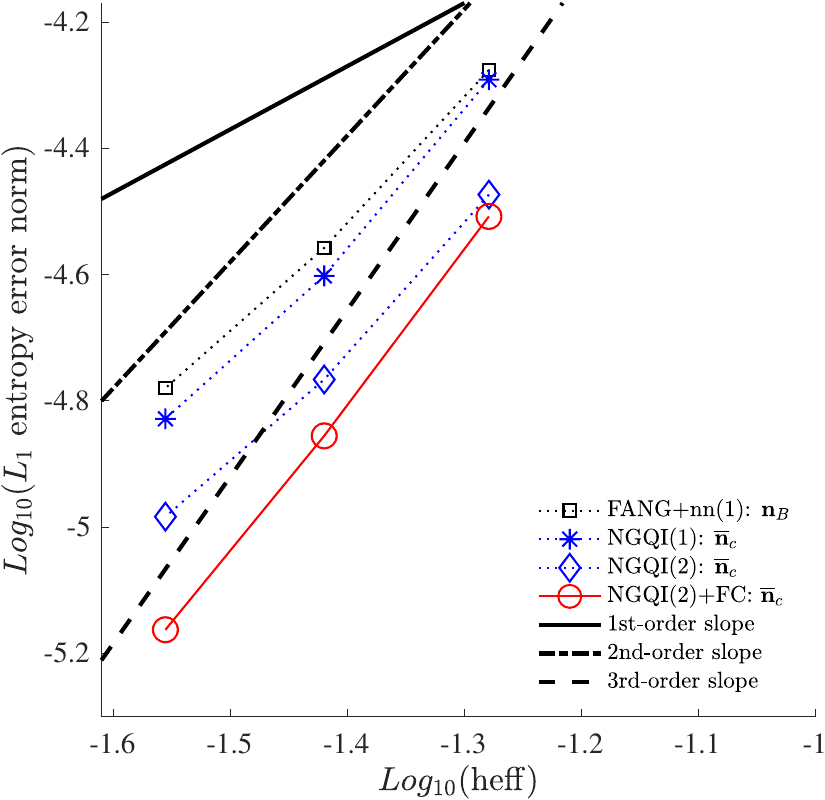}
          \caption{Entropy error convergence.}
          \label{fig:bump_err_L1_conv}
                \end{subfigure}
      \hfill 
        \caption{A subsonic flow over a smooth bump: the coarsest grid and entropy error convergence.
        }
          \label{fig:subsonic_bump}
\end{center} 
\end{figure}

To investigate the impact of the {\color{black} NGQI} scheme and the third-order scheme on a more realistic problem, we consider an $M_\infty = 0.3$ subsonic flow over a smooth bump through a duct, again with relatively coarse grids. A smooth bump defined by 
\begin{eqnarray}
z_b = 0.05 \sin^4 \left[  \frac{     \pi  ( 10 x - 3 ) }{ 9}   \right] , \quad  x \in [0.3, 1.2], 
\end{eqnarray}
is placed at the bottom of a rectangular domain $(x,y,z) = [-0.25, 1.75] \times [-0.5, 0.5 ] \times [ 0, 1] $. We solve the steady Euler equations over three levels of irregular tetrahedral grids with 9,000, 23,760, and 59,040 tetrahedra (see Figure \ref{fig:bump_grid} for the coarsest grid) by reducing the residuals by five orders of magnitude. The implicit defect-correction solver was used to solve the residual system with 10 Gauss-Seidel relaxations per iteration and the CFL number increased from 10 to 100 over the first 200 iterations. All boundary conditions are imposed weakly through the right state: the free stream condition at $x=-0.25$, the back pressure condition at $x=1.75$, and the slip-wall condition applied at other boundaries with ${\bf n}_B$ for the baseline scheme and $\hat{\bf N}_B$ for others. At a boundary node $(x_i,y_i,z_i)$ over the bump, the exact surface normal vector is given by $\hat{\bf N}_i =  {\bf N}_i / |{\bf N}_i|$, ${\bf N}_i = ( - dz_b /dx, 0, 1)$; otherwise the boundary face normal ${\bf n}_B$ is exact and we have $\hat{\bf N}_B = {\bf n}_B$. From here on, we consider a more robust version of the second-order schemes where the node-neighbor augmentation is applied to both linear and quadratic LSQ stencils as mentioned in Section \ref{LSQ_stencils}.

\begin{figure}[htbp!] 
  \begin{center}
    \centering          \begin{subfigure}[t]{0.48\textwidth}        
        \includegraphics[width=0.99\textwidth]{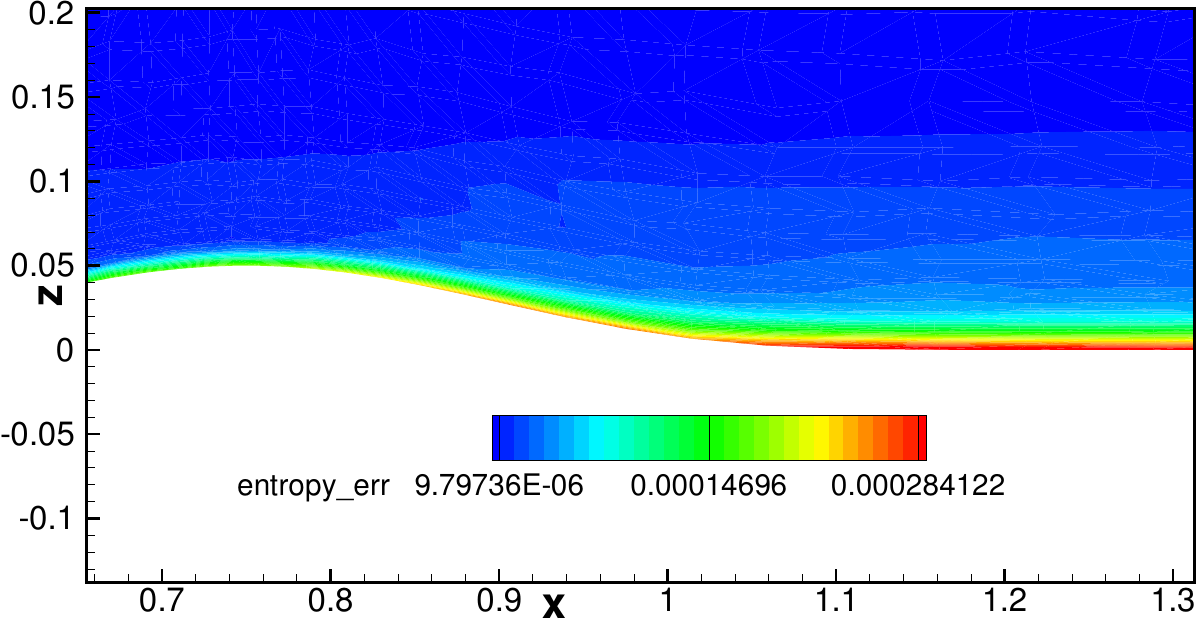}
          \caption{Entropy error contours for FANG(1): ${\bf n}_B$.}
          \label{fig:bump_grid_case00_fang_ukappa1o3_llsq_bcnorm0}
      \end{subfigure}
      \hfill
          \begin{subfigure}[t]{0.48\textwidth}        
        \includegraphics[width=0.99\textwidth]{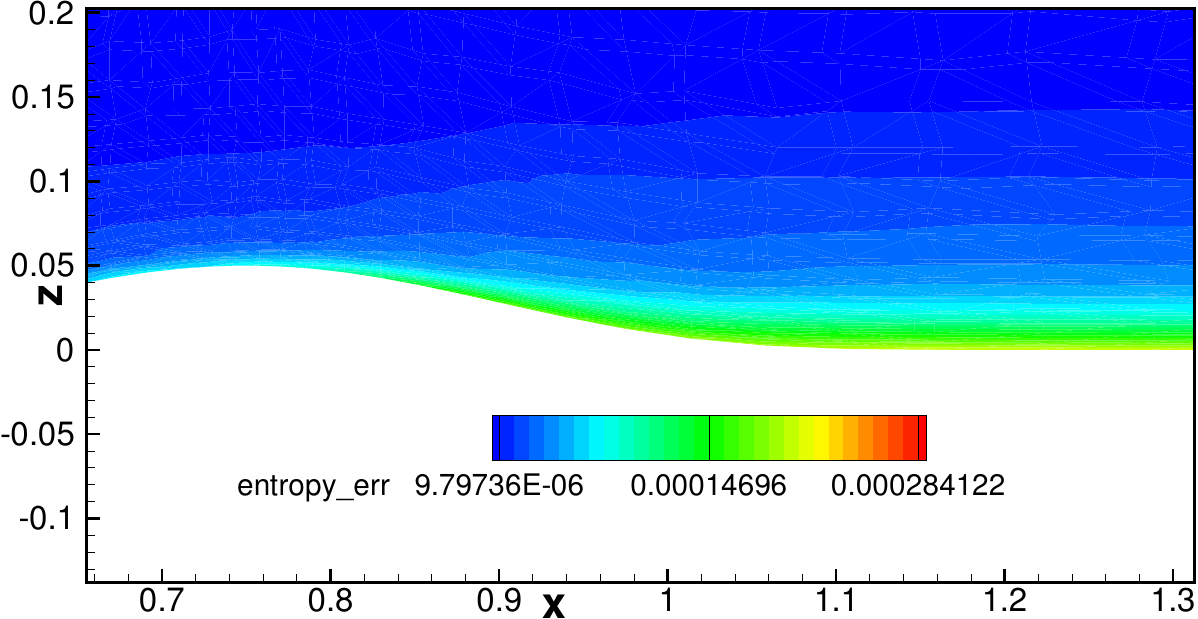}
          \caption{Entropy error contours for {\color{black} NGQI}(1): $\hat{\bf N}_B$.}
          \label{fig:bump_grid_case_effquad_llsq_bcnorm1}
      \end{subfigure}
      \hfill
          \begin{subfigure}[t]{0.48\textwidth}        
        \includegraphics[width=0.99\textwidth]{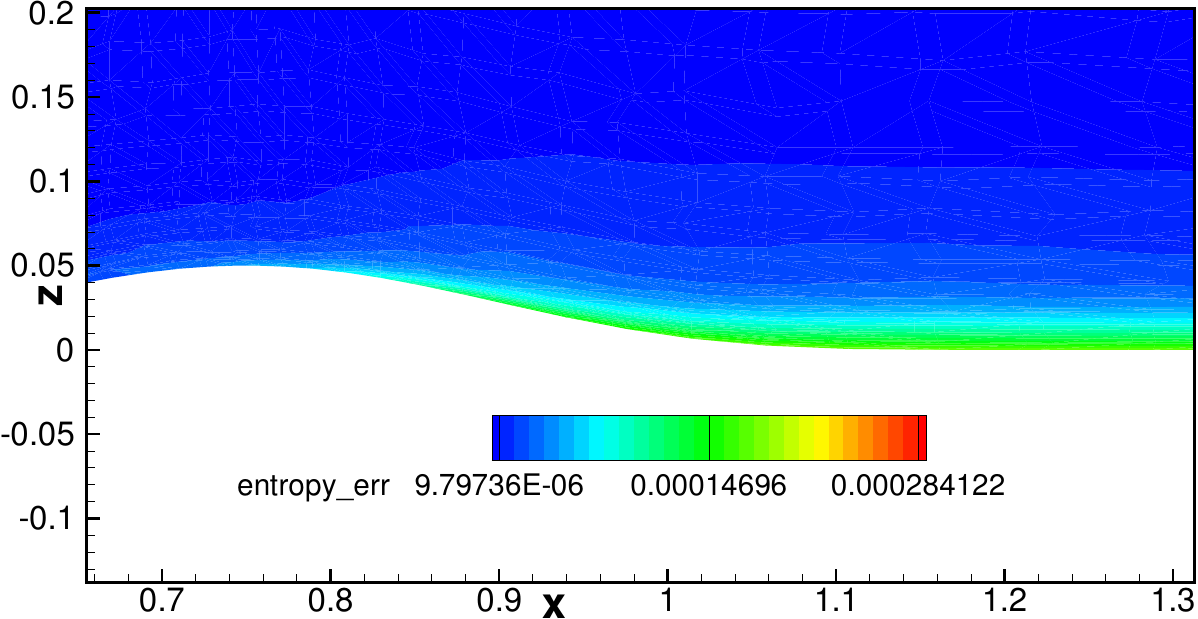}
          \caption{Entropy error contours for {\color{black} NGQI}(2): $\hat{\bf N}_B$.}
          \label{fig:bump_grid_case_effquad_qlsq_bcnorm1}
      \end{subfigure}
      \hfill
          \begin{subfigure}[t]{0.48\textwidth}        
        \includegraphics[width=0.99\textwidth]{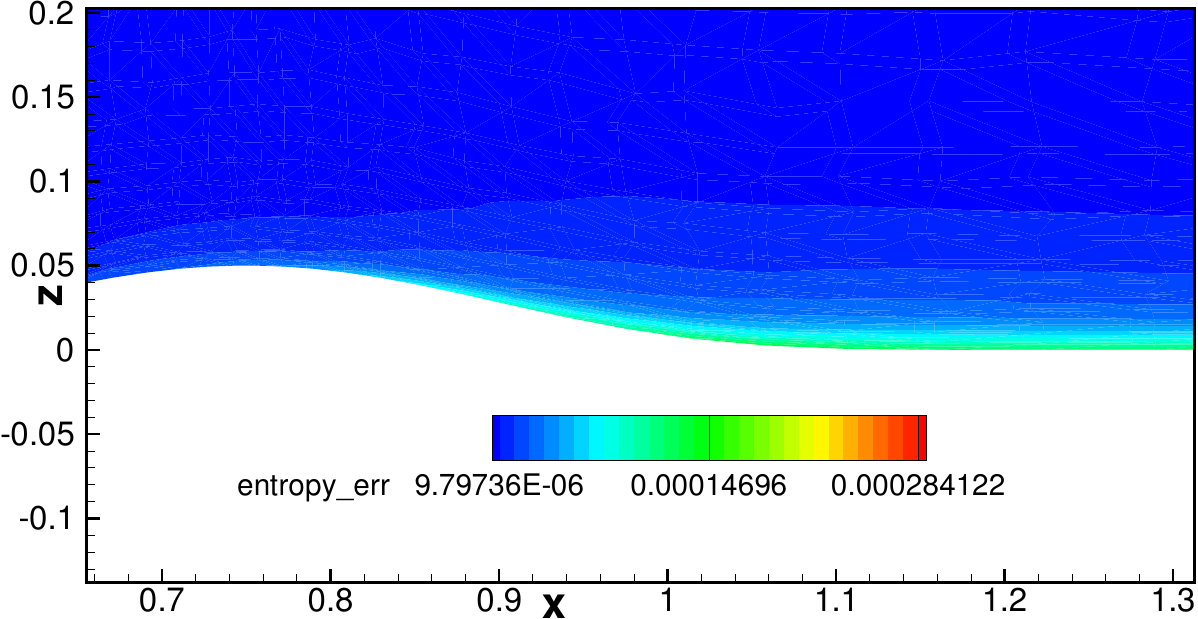}
          \caption{Entropy error contours for {\color{black} NGQI}(2)+FC: $\hat{\bf N}_B$.}
          \label{fig:bump_grid_case_effquad_cq3rd_bcnorm1}
      \end{subfigure}
      \hfill
        \caption{A subsonic flow over a smooth bump: entropy error contours.
        }
          \label{fig:subsonic_bump_entropy_err_contours_grid02p5_y0p178}
\end{center} 
\end{figure}

Figure \ref{fig:bump_err_L1_conv} shows the entropy error convergence results, where the entropy error is defined by $| \gamma p / \rho^\gamma - 1|$, which should be zero throughout the domain in such a smooth flow. As expected, the baseline scheme and {\color{black} NGQI} schemes show second-order convergence in the entropy error: {\color{black} NGQI}(2) produced a lower level of the entropy error. The third-order scheme, {\color{black} NGQI}(2)+FC, shows third-order convergence. Contours of the entropy errors are shown in Figure \ref{fig:subsonic_bump_entropy_err_contours_grid02p5_y0p178}, sampled at $y = 0.178$ on the finest grid: a lower entropy can be seen for the third-order scheme.

\begin{figure}[htbp!] 
  \begin{center}
    \centering
          \begin{subfigure}[t]{0.48\textwidth}        
        \includegraphics[width=0.99\textwidth]{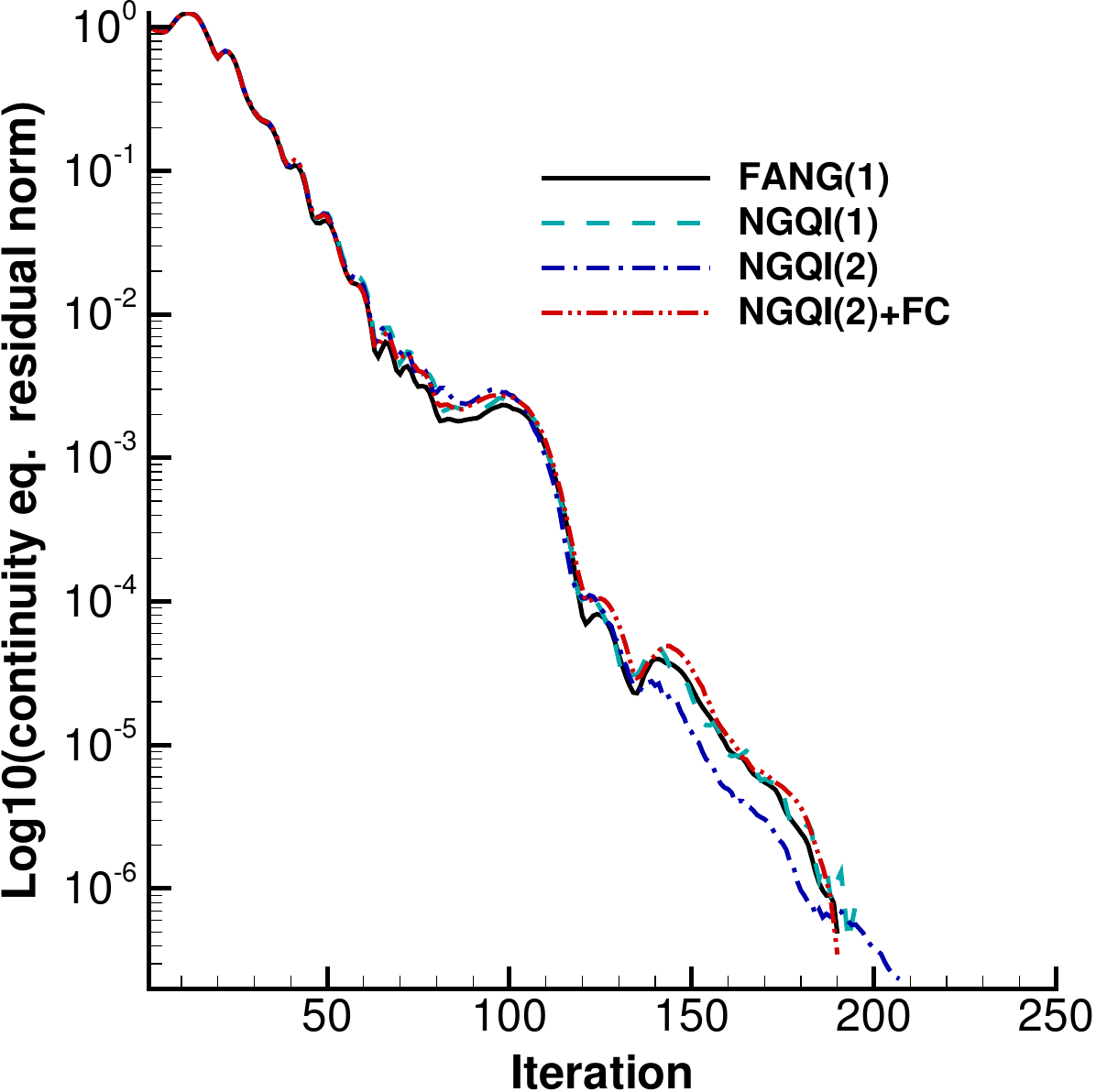}
          \caption{Residual norm versus iteration.}
          \label{fig:bump_grid_res1_vs_itr}
      \end{subfigure}
      \hfill
          \begin{subfigure}[t]{0.48\textwidth}        
        \includegraphics[width=0.99\textwidth]{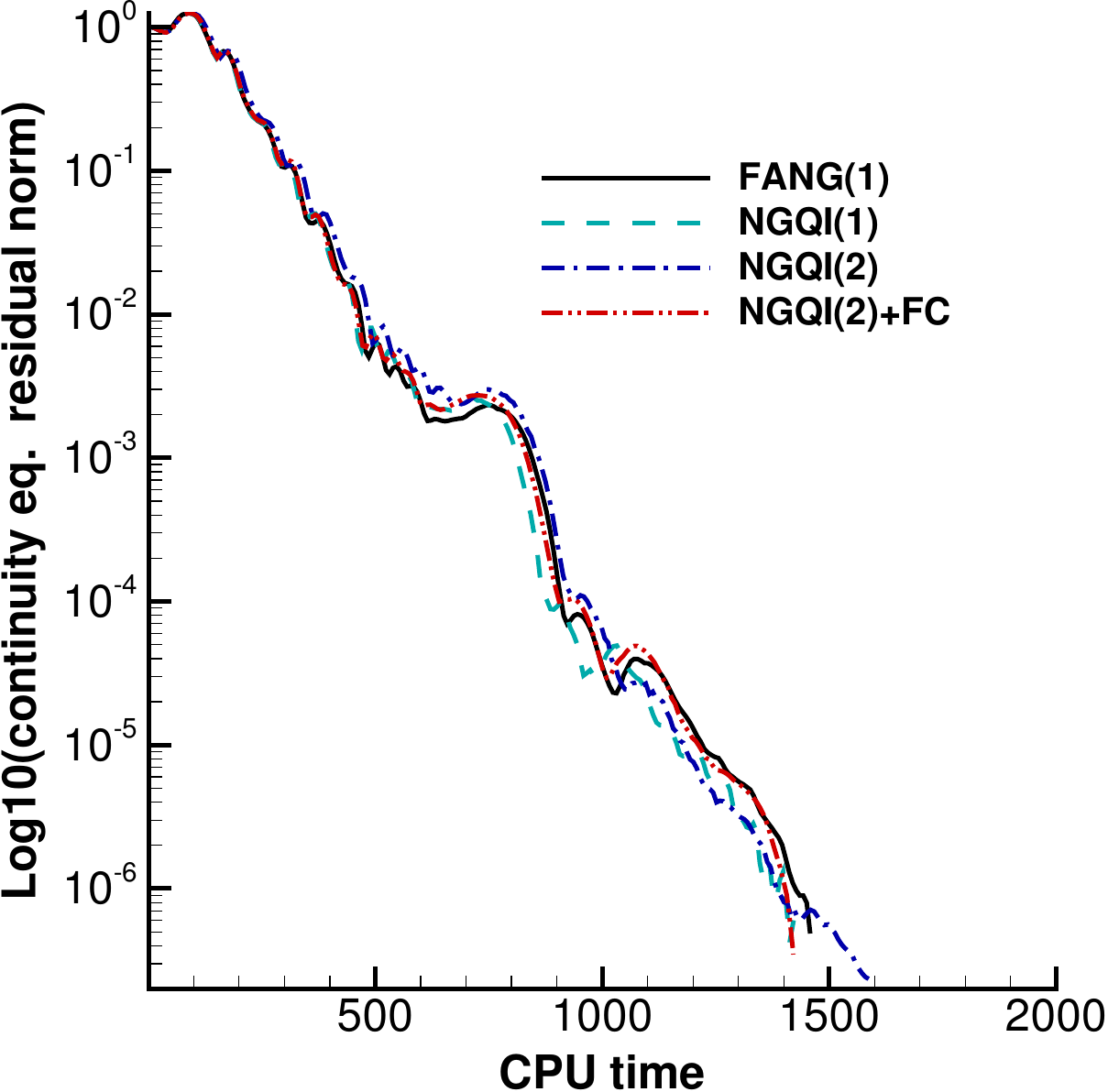}
          \caption{Residual norm versus CPU time.}
          \label{fig:bump_grid_res1_vs_wtime}
                \end{subfigure}
      \hfill 
        \caption{A subsonic flow over a smooth bump: iterative convergence.
        }
          \label{fig:subsonic_bump_itr_conv}
\end{center} 
\end{figure}

Figure \ref{fig:subsonic_bump_itr_conv} compares the iterative residual convergence for the continuity equation. As can be seen, the solver converged nearly at the same iteration for all the schemes (see Figure \ref{fig:bump_grid_res1_vs_itr}) as well as in CPU time  (see Figure \ref{fig:bump_grid_res1_vs_wtime}), indicating that the extra cost for the quadratic {\color{black} interpolation} and the flux correction is not significant, at least for the problem tested. For a steady problem, the total cost depends not only on the per-iteration cost such as the cost of residual evaluations but also the total number of iterations. It is encouraging that the third-order scheme does not require significantly more iterations to achieve superior accuracy. In the next test problem, we will compare the per-iteration cost.

\subsection{Inviscid Vortex Transport}

\begin{figure}[htbp!]
\begin{center}
\begin{minipage}[b]{0.99\textwidth}
\begin{center}
        \includegraphics[width=0.99\textwidth]{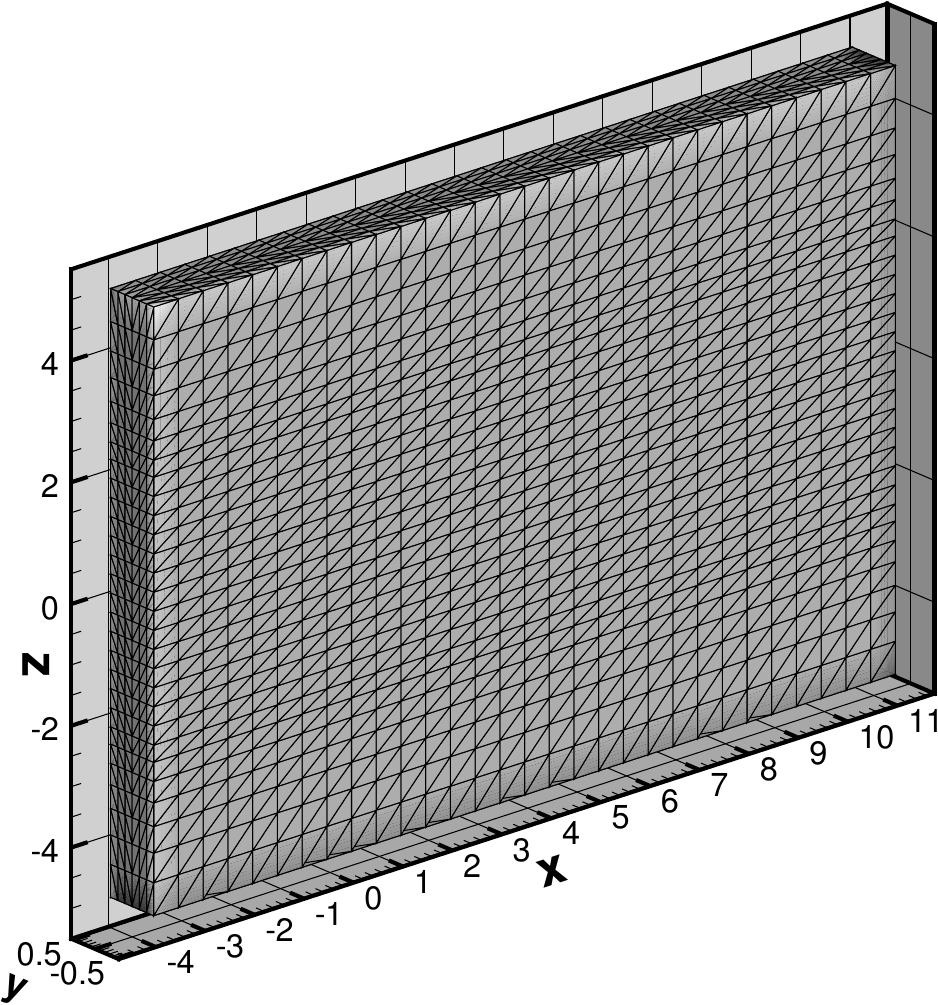}
  \caption[]{Inviscid vortex transport: a tetrahedral grid.}
   \label{fig:inv_vortex_grid}
\end{center}
\end{minipage}
\end{center}
\end{figure}

Finally, we consider the inviscid vortex transport problem, considered in Section \ref{inviscid_vortex_accuracy_verification}, to investigate the advantage of the quadratic solution {\color{black} interpolation} scheme and the third-order scheme. Again, we employ a coarse grid with 30,240 tetrahedra as shown in Figure \ref{fig:inv_vortex_grid}. Here, we will focus on the quality of the solution at a final time and the cost comparison. As in Section \ref{inviscid_vortex_accuracy_verification}, we integrate in time by using the third-order explicit SSP Runge-Kutta scheme \cite{SSP:SIAMReview2001} with a global time step $\Delta t$ = 1.3e-3 for 40,000 time steps, where the vortex initially centered at the origin will move to the point $(x,y,z) = (5.2, 0, 0)$. The global time step is chosen to be small enough (CFL $\approx$ 0.05) to focus on the comparison of the spatial discretizations. It is also relevant to potential future applications to scale-resolving simulations, in which the time step is determined by accuracy requirement rather than stability and it can be significantly smaller than the maximum allowable time step for stability. Theoretically, the proposed third-order spatial discretization might require a smaller time step than a second-order discretization for stability, but for this problem, it has been found to be stable with CFL close to 1.  Note also that a third-order explicit Runge-Kutta scheme is a suitable choice even for second-order schemes especially with a low-dissipation flux as it is known to be stable even for a second-order scheme with zero dissipation \cite{Lembert1991}. In this problem, we will focus on the same four schemes as in the previous test case, but employ the low-dissipation Roe flux \cite{nishikawa_liu_aiaa2018-4166}.

\begin{figure}[htbp!] 
  \begin{center}
    \centering
          \begin{subfigure}[t]{0.48\textwidth}        
        \includegraphics[width=0.99\textwidth]{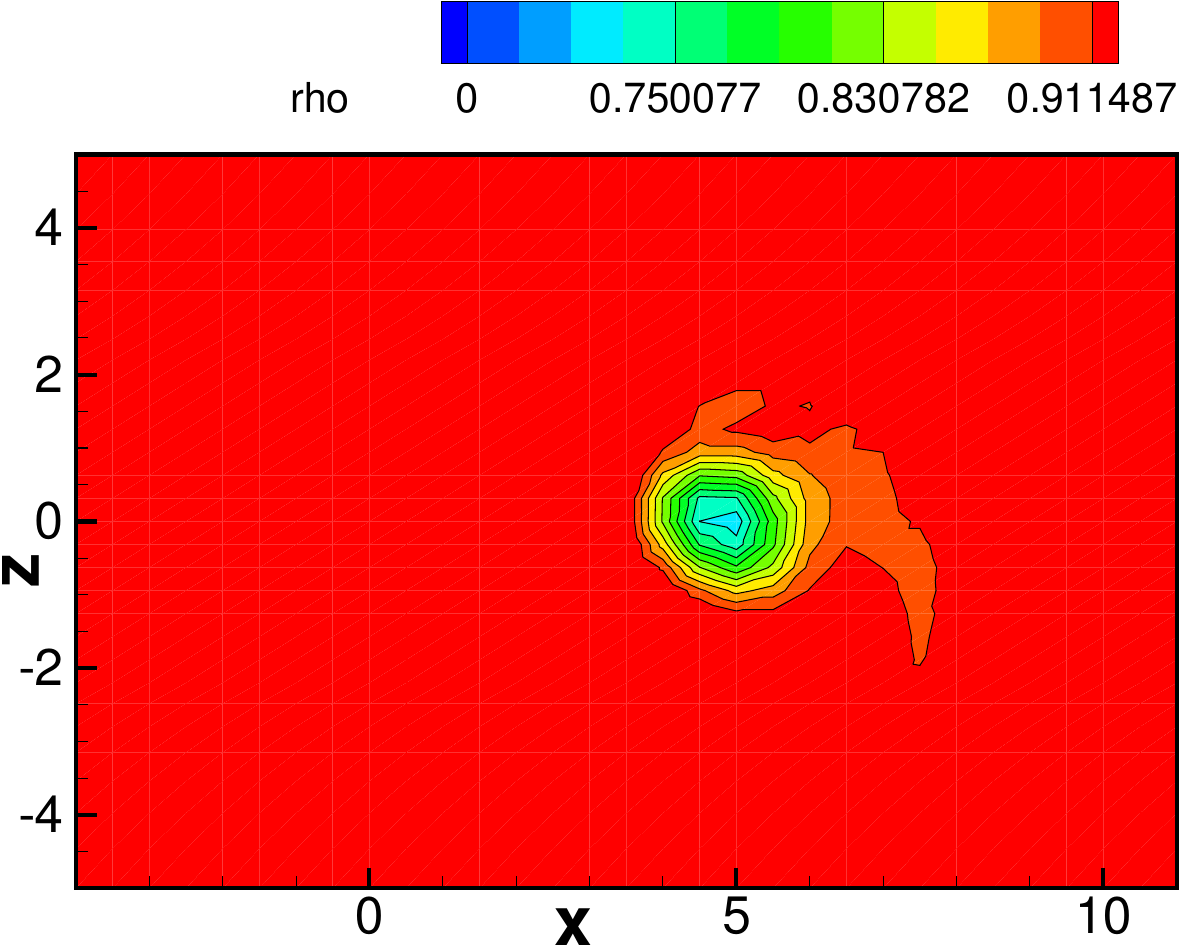}
          \caption{FANG(1).}
          \label{fig:inv_vortex_40000_fang}
      \end{subfigure}
      \hfill
          \begin{subfigure}[t]{0.48\textwidth}        
        \includegraphics[width=0.99\textwidth]{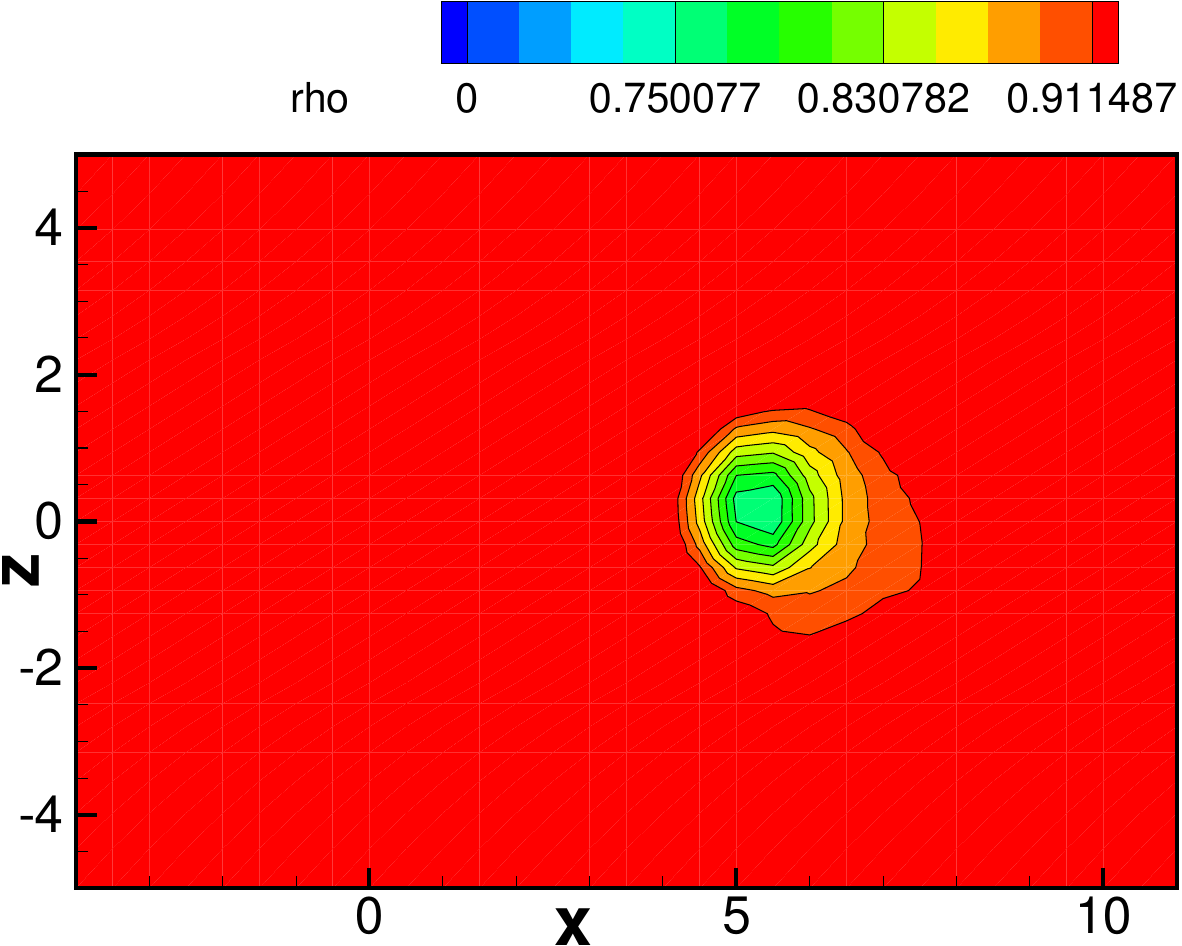}          \caption{{\color{black} NGQI}(1).}
          \label{fig:inv_vortex_40000_eqr_llsq}
                \end{subfigure}
      \hfill
          \begin{subfigure}[t]{0.48\textwidth}        
        \includegraphics[width=0.99\textwidth]{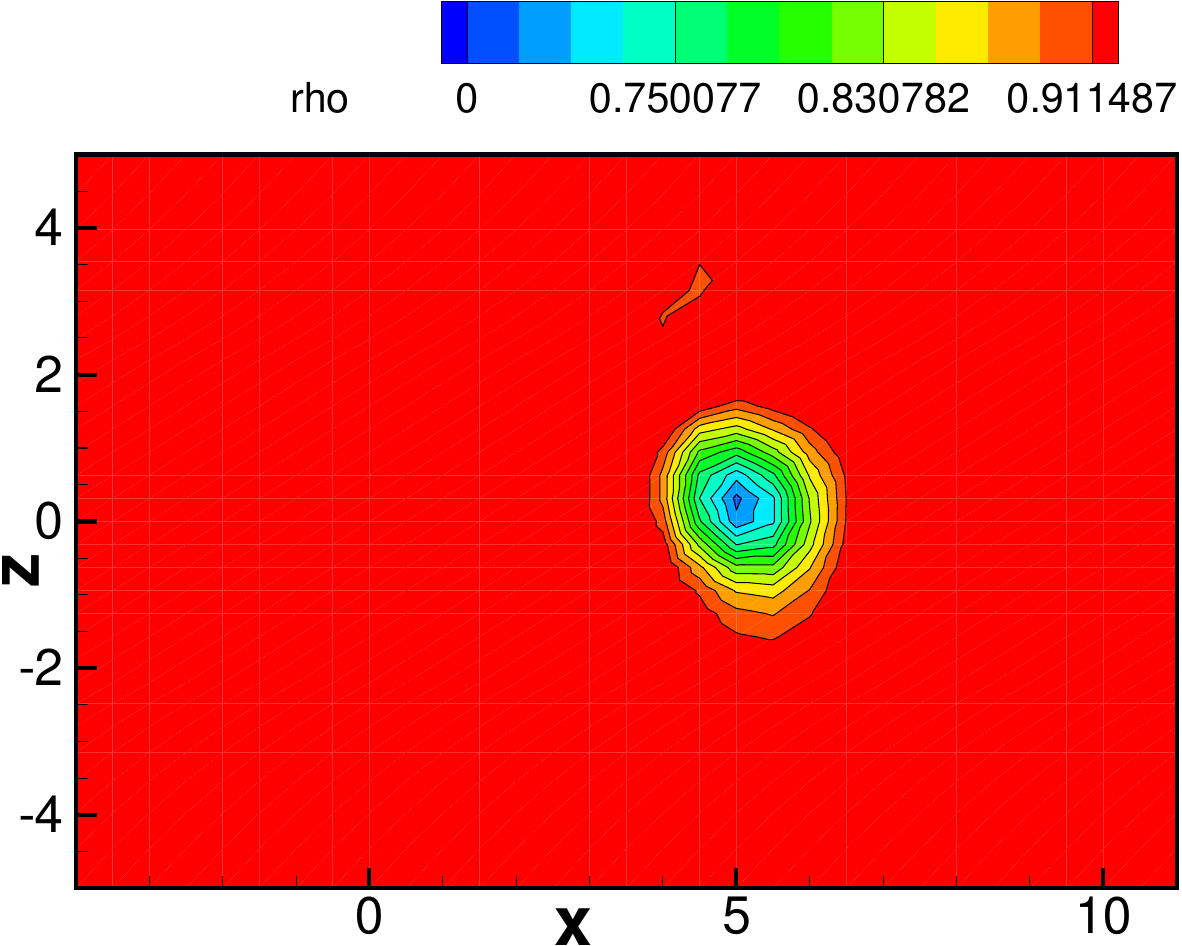}    \caption{{\color{black} NGQI}(2).}
          \label{fig:inv_vortex_40000_eqr_qlsq}
      \end{subfigure}
      \hfill
          \begin{subfigure}[t]{0.48\textwidth}        
        \includegraphics[width=0.99\textwidth]{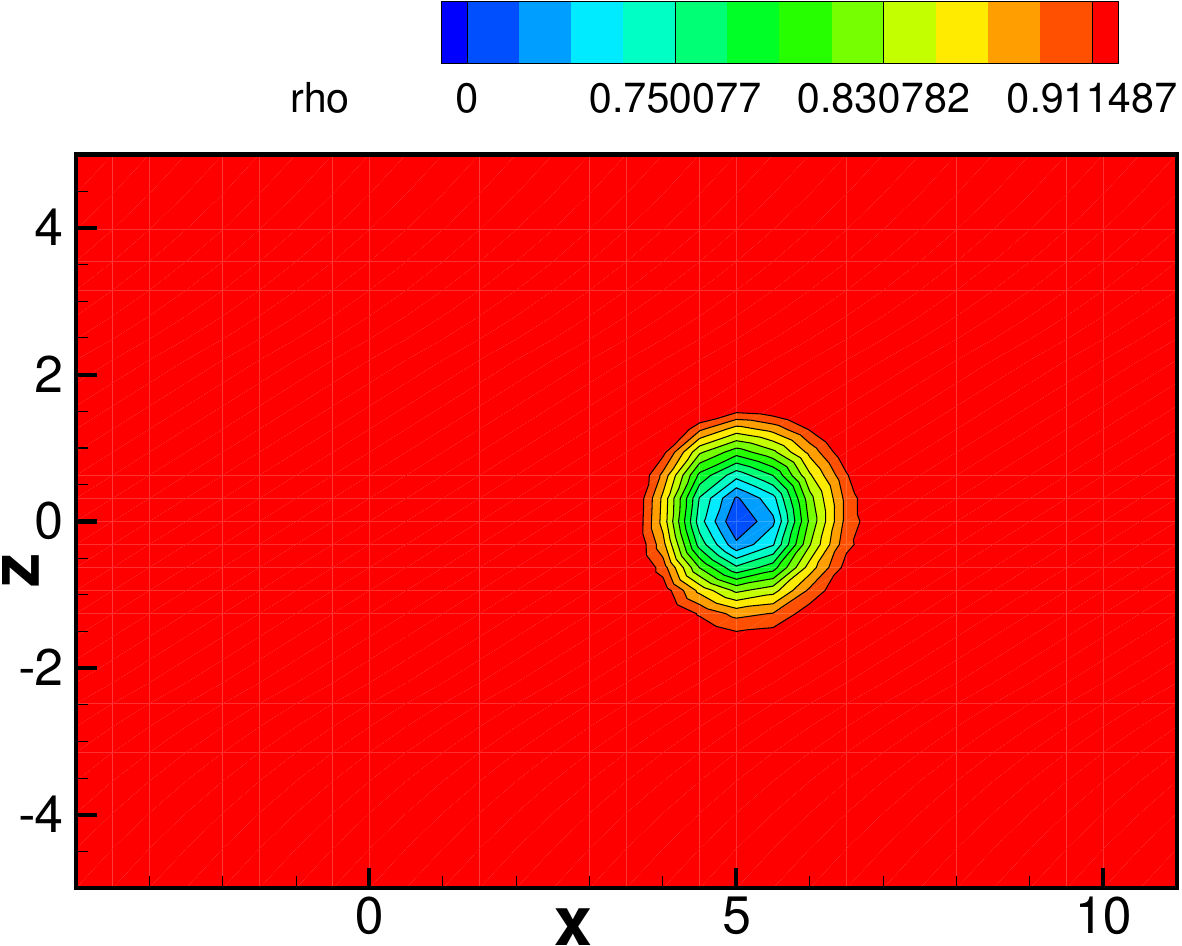}    \caption{{\color{black} NGQI}(2)+FC.}
          \label{fig:inv_vortex_40000_3rd}
      \end{subfigure}
      \hfill
          \begin{subfigure}[t]{0.48\textwidth}        
        \includegraphics[width=0.99\textwidth]{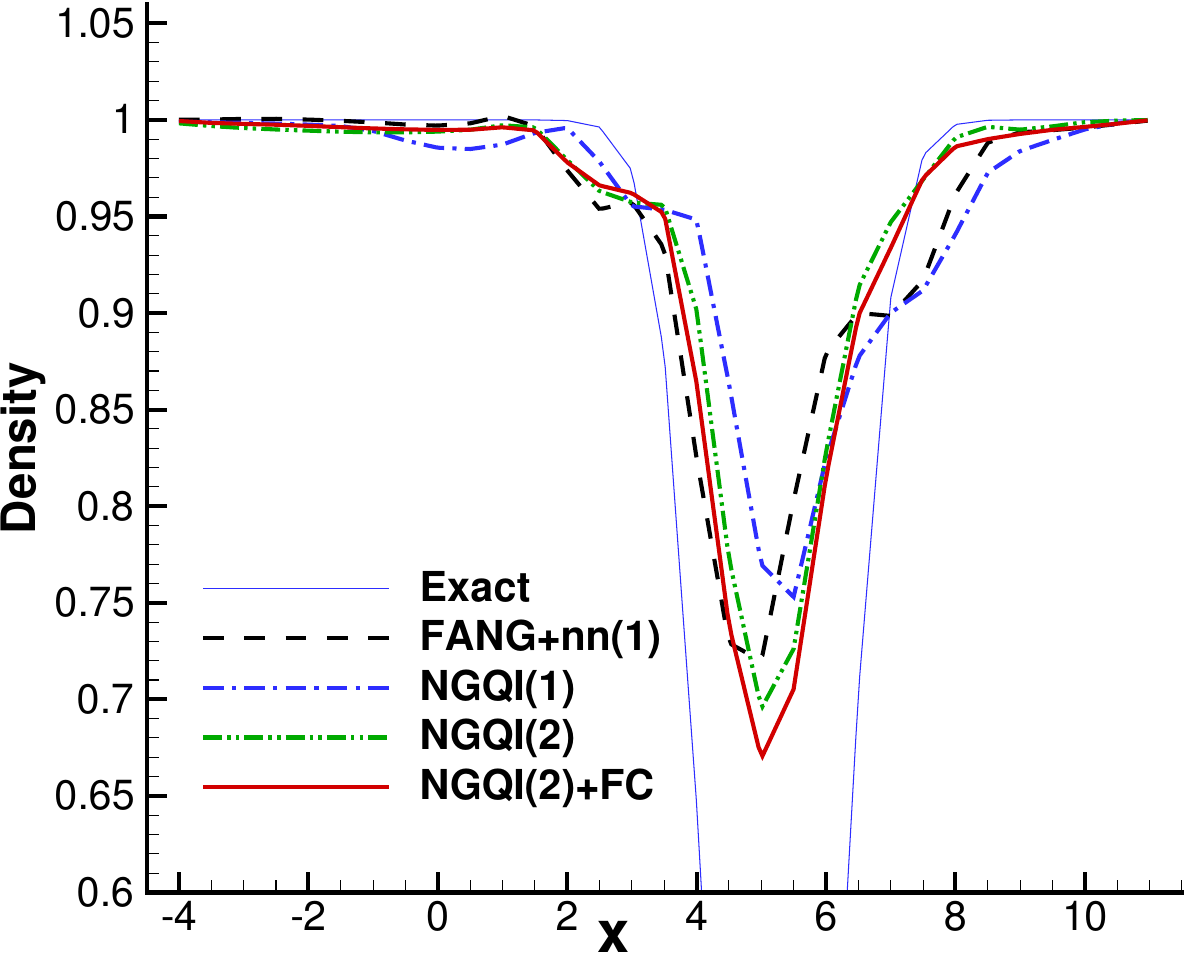}
          \caption{Section plot at $z=0.0$.}
          \label{fig:inv_vortex_40000_section}
      \end{subfigure}
      \hfill
        \caption{Inviscid vortex transport: solutions at the final time $t=5.2$ (40,000 time steps).
        }
          \label{fig:inv_vortex_40000}
\end{center} 
\end{figure}

\begin{figure}[htbp!] 
  \begin{center}
    \centering
          \begin{subfigure}[t]{0.48\textwidth}        
        \includegraphics[width=0.99\textwidth]{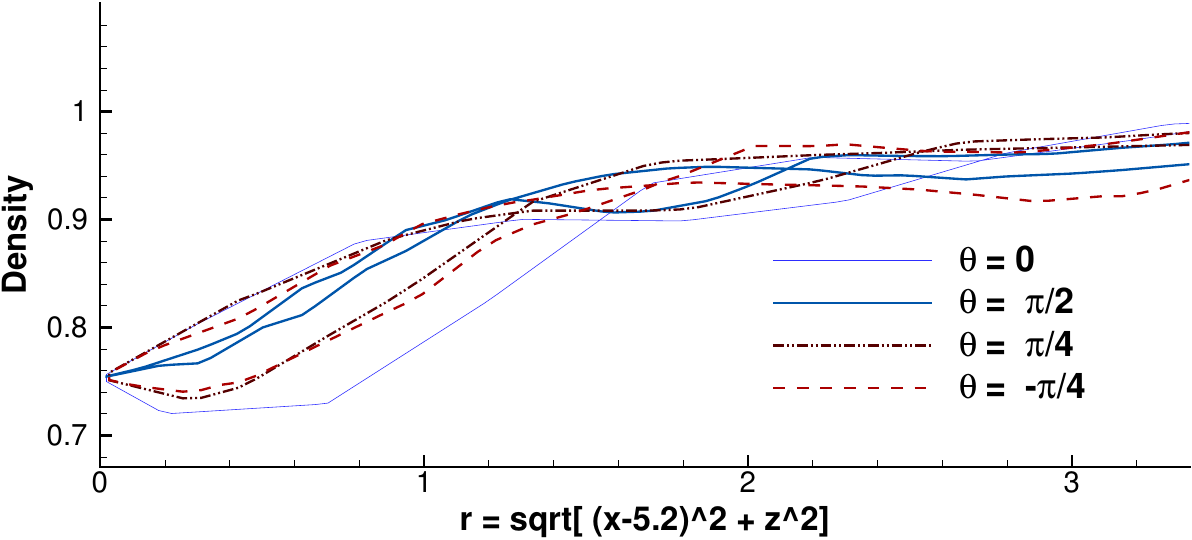}
          \caption{FANG(1).}
          \label{fig:inv_vortex_40000_fang}
      \end{subfigure}
      \hfill
          \begin{subfigure}[t]{0.48\textwidth}        
        \includegraphics[width=0.99\textwidth]{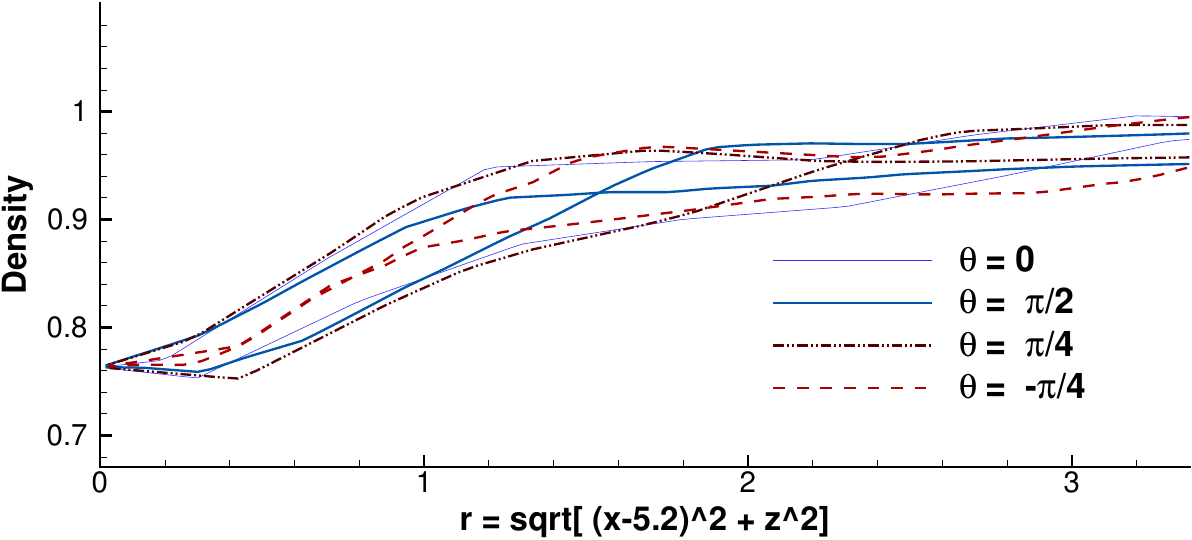}          \caption{{\color{black} NGQI}(1).}
          \label{fig:inv_vortex_40000_eqr_llsq}
                \end{subfigure}
      \hfill
          \begin{subfigure}[t]{0.48\textwidth}        
        \includegraphics[width=0.99\textwidth]{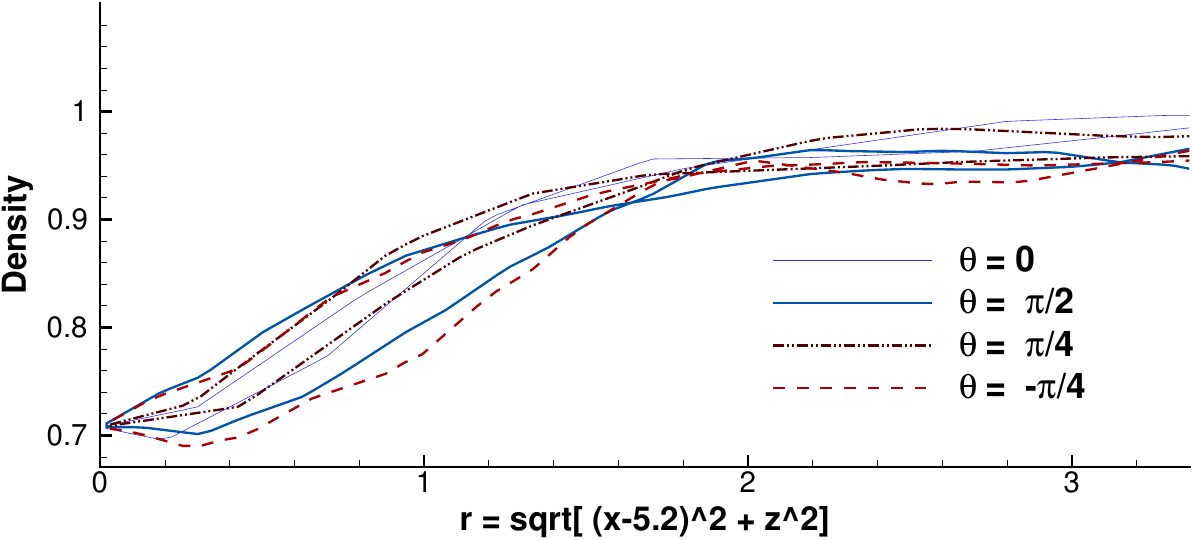} 
         \caption{{\color{black} NGQI}(2).}
          \label{fig:inv_vortex_40000_eqr_qlsq}
      \end{subfigure}
      \hfill
          \begin{subfigure}[t]{0.48\textwidth}        
        \includegraphics[width=0.99\textwidth]{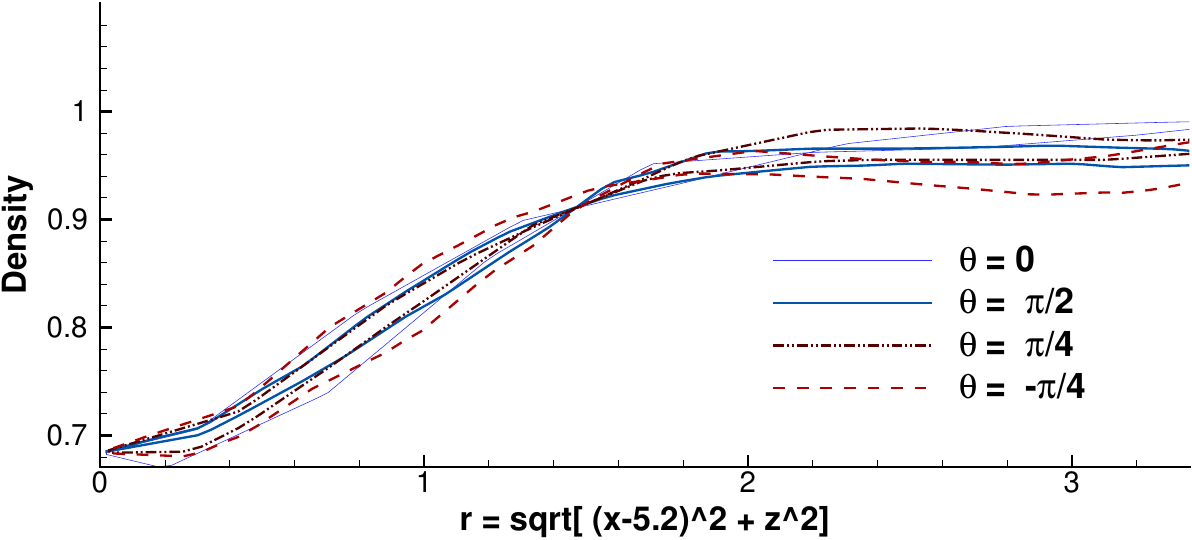}    \caption{{\color{black} NGQI}(2)+FC.}
          \label{fig:inv_vortex_40000_3rd}
      \end{subfigure}
      \hfill
        \caption{Inviscid vortex transport: section plots taken along lines at the angles $0$, $\pi/2$, $\pi/4$, and $-\pi/4$ measured from the positive $x$-axis around the center of the vortex at the final time $(5.2,0)$ (40,000 time steps).
        }
          \label{fig:inv_vortex_40000_sections}
\end{center} 
\end{figure}

Compared in Figure \ref{fig:inv_vortex_40000} are the density contour plots at 40,000 time steps ($t=52.0$). As can be seen, although all solutions are dissipative on such a coarse grid, the peak is dissipated more strongly in the baseline scheme and the {\color{black} NGQI}(1) scheme while it is slightly better preserved by the {\color{black} NGQI}(2) scheme. For the third-order scheme with the mass matrix inversion applied, {\color{black} NGQI}(2)+FC, it can be seen in Figure \ref{fig:inv_vortex_40000_3rd} that it much better preserves the peak as well as the circular shape of the solution as expected from a lower dispersive error typical of third-order schemes. Figure \ref{fig:inv_vortex_40000_section} shows the section plot taken at $z=0$. It can be seen that the peak locations of solutions obtained with the baseline scheme and the {\color{black} NGQI}(1) scheme are not accurately predicted compared with the other two schemes, indicating lower dispersion of the {\color{black} NGQI}(2) scheme and the third-order scheme. {\color{black} Figure \ref{fig:inv_vortex_40000_sections} shows the density profiles taken along lines at the angles $0$, $\pi/2$, $\pi/4$, and $-\pi/4$ measured from the positive $x$-axis around the center of the vortex at the final time $(5.2,0)$. These results show superior symmetry of the third-order solution as the density profiles are better aligned one another than those obtained with other schemes exhibiting larger variation. }

Finally, in Table \ref{Table:vortexCPUtime}, we compare the CPU time required to run the simulation with each scheme. Since all the cases were run for the same 40,000 time steps by the same time integration scheme, it basically compares the per-time-step spatial residual cost. Values indicated in the parenthesis are the relative CPU time with respect to that for the baseline scheme, FANG(1). As can be seen, {\color{black} NGQI}(1) takes almost the same CPU time as the baseline scheme while {\color{black} NGQI}(2) required about 0.8\% more CPU time. The third-order scheme, {\color{black} NGQI}(2)+FC, took around 9.6\% more CPU time. It is a bit slower than {\color{black} NGQI}(2), but the increase in cost is not significant. 
These results indicate that the new {\color{black} interpolation} scheme and the third-order scheme are efficient and can be nearly as fast as the baseline second-order scheme. 


\clearpage

\begin{table}[htbp!]
\ra{1.5}
\begin{center}
\caption{Inviscid vortex transport: total CPU time comparison.}
{
\begin{tabular}{lcccccl}\hline\hline 
\multicolumn{1}{l}{ }                                               &
\multicolumn{1}{l}{FANG(1)}                                               &
\multicolumn{1}{l}{{\color{black} NGQI}(1)}                                               &
\multicolumn{1}{l}{{\color{black} NGQI}(2)}                                             &
\multicolumn{1}{l}{{\color{black} NGQI}(2)+FC}  
  \\ \hline  
Total CPU time (sec)  &  65014 & 64938 (0.999) & 65517 (1.008) & 71276 (1.096) \\ 
   \hline  \hline
\end{tabular}
\label{Table:vortexCPUtime}
}
\end{center}
\end{table}




\section{Concluding Remarks}
\label{conclusions}

In this paper, we have developed an economical third-order cell-centered finite-volume discretization method for tetrahedral grids based on an efficient quadratic interpolation formula utilizing gradients stored at nodes. The method is unique in that it approximates {\color{black} the integral form} of a conservation law with point-valued solutions stored at the centroids, not cell averages, using an efficient one-point flux quadrature formula and a highly efficient quadratic solution interpolation scheme with no geometric moments involved. Unlike conventional third-order cell-centered finite-volume schemes, the proposed method does not require computations and storage of second derivatives of solution variables. As a result, the method can be readily implemented in an existing second-order cell-centered finite-volume code, provided gradients are computed at nodes from the solutions at cells. Superior accuracy has been demonstrated for simple but illustrative test cases in three dimensions, showing lower levels of spurious entropy errors and low-dissipation and dispersion errors for the inviscid vortex transport problem on tetrahedral grids. The proposed quadratic {\color{black} interpolation} scheme and the third-order discretization are highly efficient. For the problems considered, we have shown that these new algorithms took at most 9.6\% extra CPU time compared with a baseline second-order method.   

This paper has introduced the basic key idea for constructing an economical third-order cell-centered discretization method on tetrahedral grids. There are various possibilities for future developments. A nonoscillatory version should be developed for shock capturing. To preserve third-order accuracy in smooth regions, we would have to employ high-order limiters such as those recently developed in Ref.~\cite{Nishikawa:SciTech2022}. Extensions to viscous terms should also be investigated. First, the proposed algorithms are directly applicable to the hyperbolic Navier-Stokes formulations \cite{liu_nishikawa_aiaa2016-3969,nishikawa_hyperbolic_ns:AIAA2015,NakashimaWatanabeNishikawa_AIAA2016-1101,LiLouNishikawaLuo_HNSrDG:JCP2021}; then, third-order accuracy would be achieved in the solution variables as well as in their gradients (e.g., viscous stresses and heat fluxes) on irregular grids. For conventional second-derivative viscous terms, however, it remains to be investigated whether third-order accuracy can be achieved with the quadratic {\color{black} interpolation}; other studies indicate that a quadratic {\color{black} interpolation} is not sufficient to achieve third-order accuracy for the viscous terms \cite{SetzweinEssGerlinger:JCP2021} and a cubic {\color{black} interpolation} is required to achieve third-order accuracy for the viscous terms \cite{pincock_katz:JSC2014_DOI}. {\color{black} Extensions to mixed-element grids are also important from a practical point of view and currently underway. In principle, the efficient quadratic {\color{black} interpolation} scheme can be extended to cells other than tetrahedra; but it might be necessary to store interpolation coefficients in each cell for computing $\overline{\bf g}_j$ in Eq.~(\ref{efficient_quadratic_reconstruction_final_wL_xT}).} Another interesting future development would be the use of implicit gradient methods \cite{WangRenPanLi:JCP2017,nishikawa_aiaa2019-1155,Nishikawa_IEBG:Aviation2020}, which could reduce the cost of gradient computation as well as improve nonlinear solver stability. Finally, the method should be demonstrated for realistic three-dimensional turbulent-flow problems by a practical parallel code. 

\addcontentsline{toc}{section}{Acknowledgments}
\section*{Acknowledgments}

The authors gratefully acknowledge support by the Hypersonic Technology Project, through the Hypersonic Airbreathing Propulsion Branch of the NASA Langley Research Center. The first author was funded under Contract No. 80LARC17C0004. The first author would like to thank Boris Diskin (National Institute of Aerospace) for a discussion on accuracy of slip boundary conditions. {\color{black} We would like to thank reviewers for their constructive comments.}

\bibliography{../../../..//bibtex_nishikawa_database}
\bibliographystyle{unsrt}




\end{document}